\documentclass[12pt]{article}
\usepackage{latexsym, amssymb}

\usepackage{mathbbol}

\makeatletter
\def\eqnarray{\stepcounter{equation}\let\@currentlabel=\theequation
\global\@eqnswtrue
\tabskip\@centering\let\\=\@eqncr
$$\halign to \displaywidth\bgroup\hfil\global\@eqcnt\z@
  $\displaystyle\tabskip\z@{##}$&\global\@eqcnt\@ne
  \hfil$\displaystyle{{}##{}}$\hfil
  &\global\@eqcnt\tw@ $\displaystyle{##}$\hfil
  \tabskip\@centering&\llap{##}\tabskip\z@\cr}

\def\endeqnarray{\@@eqncr\egroup
      \global\advance\c@equation\m@ne$$\global\@ignoretrue}

\def\@yeqncr{\@ifnextchar [{\@xeqncr}{\@xeqncr[5pt]}}
\makeatother

\begin{document}
\bibliographystyle{tom}

\newtheorem{lemma}{Lemma}[section]
\newtheorem{thm}[lemma]{Theorem}
\newtheorem{cor}[lemma]{Corollary}
\newtheorem{voorb}[lemma]{Example}
\newtheorem{rem}[lemma]{Remark}
\newtheorem{prop}[lemma]{Proposition}
\newtheorem{ddefinition}[lemma]{Definition}
\newtheorem{stat}[lemma]{{\hspace{-5pt}}}

\newenvironment{remarkn}{\begin{rem} \rm}{\end{rem}}
\newenvironment{exam}{\begin{voorb} \rm}{\end{voorb}}
\newenvironment{definition}{\begin{ddefinition} \rm}{\end{ddefinition}}

\newcounter{teller}
\renewcommand{\theteller}{(\alph{teller})}
\newenvironment{tabel}{\begin{list}%
{\rm  (\alph{teller})\hfill}{\usecounter{teller} \leftmargin=1.1cm
\labelwidth=1.1cm \labelsep=0cm \parsep=0cm}
                      }{\end{list}}

\newcounter{tellerr}
\renewcommand{\thetellerr}{(\roman{tellerr})}
\newenvironment{tabeleq}{\begin{list}%
{\rm  (\roman{tellerr})\hfill}{\usecounter{tellerr} \leftmargin=1.1cm
\labelwidth=1.1cm \labelsep=0cm \parsep=0cm}
                         }{\end{list}}

\newcounter{tellerrr}
\renewcommand{\thetellerrr}{(\Roman{tellerrr})}
\newenvironment{tabelR}{\begin{list}%
{\rm  (\Roman{tellerrr})\hfill}{\usecounter{tellerrr} \leftmargin=1.1cm
\labelwidth=1.1cm \labelsep=0cm \parsep=0cm}
                         }{\end{list}}

\newcounter{proofstep}
\newcommand{\nextstep}{\refstepcounter{proofstep}\ruimte \par 
          \noindent{\bf Step \theproofstep} \hspace{5pt}}
\newcommand{\firststep}{\setcounter{proofstep}{0}\nextstep}

\newcommand{\Ni}{{\bf N}}
\newcommand{\Qi}{{\bf Q}}
\newcommand{\Ri}{{\bf R}}
\newcommand{\Ci}{{\bf C}}
\newcommand{\Ti}{{\bf T}}
\newcommand{\Zi}{{\bf Z}}
\newcommand{\Fi}{{\bf F}}

\newcommand{\proof}{\mbox{\bf Proof} \hspace{5pt}} 
\newcommand{\remark}{\mbox{\bf Remark} \hspace{5pt}}
\newcommand{\ruimte}{\vskip10.0pt plus 4.0pt minus 6.0pt}

\newcommand{\simh}{{\stackrel{{\rm cap}}{\sim}}}
\newcommand{\ad}{{\mathop{\rm ad}}}
\newcommand{\Ad}{{\mathop{\rm Ad}}}
\newcommand{\Aut}{\mathop{\rm Aut}}
\newcommand{\arccot}{\mathop{\rm arccot}}
\newcommand{\capp}{{\mathop{\rm cap}}}
\newcommand{\rcapp}{{\mathop{\rm rcap}}}
\newcommand{\diam}{\mathop{\rm diam}}
\newcommand{\divv}{\mathop{\rm div}}
\newcommand{\codim}{\mathop{\rm codim}}
\newcommand{\RRe}{\mathop{\rm Re}}
\newcommand{\IIm}{\mathop{\rm Im}}
\newcommand{\Tr}{{\mathop{\rm Tr \,}}}
\newcommand{\Vol}{{\mathop{\rm Vol}}}
\newcommand{\card}{{\mathop{\rm card}}}
\newcommand{\supp}{\mathop{\rm supp}}
\newcommand{\sgn}{\mathop{\rm sgn}}
\newcommand{\essinf}{\mathop{\rm ess\,inf}}
\newcommand{\esssup}{\mathop{\rm ess\,sup}}
\newcommand{\Int}{\mathop{\rm Int}}
\newcommand{\lcm}{\mathop{\rm lcm}}
\newcommand{\loc}{{\rm loc}}

\newcommand{\mod}{\mathop{\rm mod}}
\newcommand{\spann}{\mathop{\rm span}}
\newcommand{\one}{\mathbb{1}}

\hyphenation{groups}
\hyphenation{unitary}

\newcommand{\tfrac}[2]{{\textstyle \frac{#1}{#2}}}

\newcommand{\ca}{{\cal A}}
\newcommand{\cb}{{\cal B}}
\newcommand{\cc}{{\cal C}}
\newcommand{\cd}{{\cal D}}
\newcommand{\ce}{{\cal E}}
\newcommand{\cf}{{\cal F}}
\newcommand{\ch}{{\cal H}}
\newcommand{\ci}{{\cal I}}
\newcommand{\ck}{{\cal K}}
\newcommand{\cl}{{\cal L}}
\newcommand{\cm}{{\cal M}}
\newcommand{\co}{{\cal O}}
\newcommand{\cs}{{\cal S}}
\newcommand{\ct}{{\cal T}}
\newcommand{\cx}{{\cal X}}
\newcommand{\cy}{{\cal Y}}
\newcommand{\cz}{{\cal Z}}

\newlength{\hightcharacter}
\newlength{\widthcharacter}
\newcommand{\covsup}[1]{\settowidth{\widthcharacter}{$#1$}\addtolength{\widthcharacter}{-0.15em}\settoheight{\hightcharacter}{$#1$}\addtolength{\hightcharacter}{0.1ex}#1\raisebox{\hightcharacter}[0pt][0pt]{\makebox[0pt]{\hspace{-\widthcharacter}$\scriptstyle\circ$}}}
\newcommand{\cov}[1]{\settowidth{\widthcharacter}{$#1$}\addtolength{\widthcharacter}{-0.15em}\settoheight{\hightcharacter}{$#1$}\addtolength{\hightcharacter}{0.1ex}#1\raisebox{\hightcharacter}{\makebox[0pt]{\hspace{-\widthcharacter}$\scriptstyle\circ$}}}
\newcommand{\scov}[1]{\settowidth{\widthcharacter}{$#1$}\addtolength{\widthcharacter}{-0.15em}\settoheight{\hightcharacter}{$#1$}\addtolength{\hightcharacter}{0.1ex}#1\raisebox{0.7\hightcharacter}{\makebox[0pt]{\hspace{-\widthcharacter}$\scriptstyle\circ$}}}

\thispagestyle{empty}

\vspace*{1cm}
\begin{center}
{\Large\bf Sectorial forms and degenerate differential operators} \\[5mm]

\large W. Arendt$^1$ and A.F.M. ter Elst$^2$

\end{center}

\vspace{5mm}

\begin{center}
{\bf Abstract}
\end{center}

\begin{list}{}{\leftmargin=1.8cm \rightmargin=1.8cm \listparindent=10mm 
   \parsep=0pt}
\item
If $a$ is a densely defined sectorial form in a Hilbert space which is possibly
not closable, then we associate in a natural way a holomorphic semigroup
generator with~$a$.
This allows us to remove in several theorems of semigroup theory 
the assumption that the form is closed or symmetric. 
Many examples are provided, ranging from complex sectorial differential
operators, to Dirichlet-to-Neumann operators and operators with 
Robin or Wentzell boundary conditions.
\end{list}

\vspace{4cm}
\noindent
April 2010

\vspace{5mm}

\vspace{5mm}

\vspace{15mm}

\noindent
{\bf Home institutions:}    \\[3mm]
\begin{tabular}{@{}cl@{\hspace{10mm}}cl}
1. & Institute of Applied Analysis  & 
  2. & Department of Mathematics   \\
& University of Ulm   & 
  & University of Auckland   \\
& Helmholtzstr.\ 18 & 
  & Private bag 92019  \\
& 89081 Ulm & 
  & Auckland  \\
& Germany  & 
  & New Zealand  \\[8mm]
\end{tabular}

\newpage
\setcounter{page}{1}

\section{Introduction} \label{Ssesqui1}

Form methods are most efficient to solve evolution equations in a 
Hilbert space~$H$.
The theory establishes a correspondence between closable sectorial forms 
and holomorphic semigroups on $H$ which are contractive on 
a sector (see Kato \cite{Kat1}, Tanabe \cite{Tan} and Ma--R\"ockner
\cite{MR}, for example).
The aim of this article is to extend the theory in two directions
and apply the new 
criteria to differential operators.
Our first result shows that the condition of closability
can be omitted completely.
To be more precise, consider a sesquilinear form 
\[
a \colon D(a) \times D(a) \to \Ci
\]
where $D(a)$ is a dense subspace of a Hilbert space~$H$.
The form a is called {\bf sectorial} if there exist a (closed) sector
\[
\Sigma_\theta
= \{ r \, e^{i \alpha} : r \geq 0, \; |\alpha| \leq \theta \}
\]
with $\theta \in [0,\frac{\pi}{2})$, and $\gamma \in \Ri$, such that 
$a(u) - \gamma \, \|u\|_H^2 \in \Sigma_\theta$ for all $u \in D(a)$,
where $a(u) = a(u,u)$.
We shall show that there exists an operator $A$ in $H$ such that 
for all $x,f \in H$ one has
$x \in D(A)$ and $Ax = f$ if and only if
there exist $u_1,u_2,\ldots \in D(a)$ such that $(\RRe a(u_n))_n$ is bounded, 
$\lim_{n \to \infty} u_n = x$ in $H$
and $\lim_{n \to \infty} a(u_n,v) = (f,v)_H$ for all $v\in D(a)$.
It is  part of the following theorem that $f$ is independent of the  
sequence $u_1,u_2,\ldots$.

\begin{thm}[Incomplete case] \label{tsesqui322}
The operator $A$ is well-defined and $-A$ generates a 
holomorphic $C_0$-semigroup 
on the interior of $\Sigma_{\frac{\pi}{2} - \theta}$.
\end{thm}

This is a special case of Theorem~\ref{tsesqui340} below, but we give a 
short proof already in Section~\ref{Ssesqui2}.
Recall that the form $a$ is called {\bf closable} 
if for every Cauchy sequence $u_1,u_2,\ldots$
in $D(a)$ such that $\lim_{n \to \infty} u_n = 0$
in $H$ one has $\lim_{n \to \infty} a(u_n) = 0$.
Here $D(a)$ carries the natural norm
$\|u\|_a = (\RRe a(u) + (1 - \gamma) \, \|u\|_H^2)^{1/2}$.
In Theorem~\ref{tsesqui322} we do not assume that $a$ is closable.
Nonetheless, the operator $A$ is well-defined.

For our second extension of the theory on form methods we consider the 
{\bf complete case}, where the form $a$ is defined 
on a Hilbert space $V$.
However, we do not assume that $V$ is embedded in~$H$, but merely
that there exists
a not necessarily injective operator $j$ from $V$ into~$H$.
This case is actually the first we consider in Section~\ref{Ssesqui2}.
It is used for the proof of Theorem~\ref{tsesqui322} given 
in Section~\ref{Ssesqui2}.
In Theorem~\ref{tsesqui340} we give a common extension of both
Theorem~\ref{tsesqui322} and the main theorem of Section~\ref{Ssesqui2}.
It turns out that many examples can be treated by our extended form method
and Section~\ref{Ssesqui4} is devoted to several applications.
Our most substantial results 
concern degenerate 
elliptic differential operators of second order with 
complex measurable coefficients on an open set $\Omega$ in $\Ri^d$.
If the coefficients satisfy merely a sectoriality condition 
(which can be very degenerate including the case where the 
coefficients are zero on some part of  $\Omega$),
then Theorem~\ref{tsesqui322} shows right away 
that the corresponding operator generates  
a holomorphic $C_0$-semigroup on $L_2(\Omega)$.
We prove a Davies--Gaffney type estimate
which gives us locality properties and in case of 
Neumann boundary conditions and real coefficients, 
the invariance of the constant functions.
This extends results for positive symmetric
forms on $\Ri^d$ in \cite{ERSZ1} and \cite{ERSZ2}.
We also extend the criteria for closed convex sets due to Ouhabaz
\cite{Ouh5} to our more general situation and show that the 
semigroup is submarkovian if the coefficients are real
(but possibly non-symmetric).
As a second application, we present an easy and direct 
treatment of the Dirichlet-to-Neumann operator
on a Lipschitz domain~$\Omega$.
Here it is essential to allow non injective $j \colon D(a) \to H$.
As a result, we obtain submarkovian semigroups on $L_p(\partial \Omega)$.
Most interesting are Robin boundary conditions which we consider
in Subsection~\ref{Ssesqui4.5} on an open bounded
set $\Omega$ of $\Ri^d$ with the 
$(d-1)$-dimensional Hausdorff measure on $\partial \Omega$.
Using Theorem~\ref{tsesqui322} we obtain directly a 
holomorphic semigroup on $L_2(\Omega)$.
Moreover, for every element in the domain of the 
generator there is a unique trace in $L_2(\partial \Omega, \sigma)$
realising Robin boundary conditions.
Such boundary conditions on rough domains had been considered
before by Daners \cite{Daners2} and \cite{AW2}.
We also give a new simple proof for the existence of a trace
for such general domains.
We use these results on the trace to 
consider Wentzell boundary conditions in Subsection~\ref{Ssesqui4.3}.
These boundary conditions obtained much attention recently
\cite{FGGR} \cite{VV}.
By our approach we may allow degenerate coefficients for the 
elliptic operator and the boundary condition.
Our final application in Subsection~\ref{Ssesqui4.4} concerns
multiplicative perturbation of the Laplacian.

Throughout this paper we use the notation and conventions as in 
\cite{Kat1}.
Moreover, the field is $\Ci$, except if indicated explicitly.
We will only consider univocal operators.

\section{Generation theorems for the complete case} \label{Ssesqui2}

The first step in the proof of Theorem~\ref{tsesqui322} is the following 
extension of the `French' approach to closed sectorial forms
(see Dautray--Lions \cite{DL5} Chapter~XVIIA Example~3,
Tanabe \cite{Tan} Sections~2.2 and 3.6, and 
Lions \cite{Lio2}).
It is a generation theorem for forms with a complete form domain. 
It differs from the usual well-known result for closed forms in the following point. 
We do not assume that the form domain is a subspace of the given 
Hilbert space, but that there exists a linear mapping $j$ from the form domain 
into the Hilbert space. 
Moreover, we do not assume that the mapping is injective. 
In the injective case, and also in the general case by restricting $j$ to 
the orthogonal complement of its kernel, we could reduce our 
result to the usual case. 
It seems to us simpler to give a direct proof, though, 
which is adapted from \cite{Tan}, Section~3.6, Application~2, 
treating the usual case.

Let $V$ be a normed space and $a \colon V \times V \to \Ci$ a 
sesquilinear form.
Then $a$ is continuous if and only if
there exists a $c > 0$ such that 
\begin{equation}
|a(u,v)| \leq c \, \|u\|_V \, \|v\|_V
\label{eSsesqui2.3}
\end{equation}
for all $u,v \in V$.
Let $H$ be a Hilbert space and $j\colon V \to H$ a bounded linear operator.
The form $a \colon V \times V \to \Ci$ is called 
$j$-{\bf elliptic} if there exist $\omega \in \Ri$
and $\mu > 0$ such that
\begin{equation}
\RRe a(u) + \omega \, \|j(u)\|_H^2 \geq \mu \, \|u\|_V^2
\label{etsesqui201;1}
\end{equation}
for all $u \in V$.
The form $a$ is called {\bf coercive} if (\ref{etsesqui201;1}) is valid 
with $\omega = 0$.

An operator $A \colon D(A) \to H$ with $D(A) \subset H$ is called
{\bf sectorial}  if 
there are $\gamma \in \Ri$, called a {\bf vertex},
and $\theta \in [0,\frac{\pi}{2})$, called a {\bf semi-angle}, such that 
\[
(Ax,x) - \gamma \, \|x\|_H^2
\in \Sigma_\theta
\]
for all $x \in D(A)$.
Moreover, $A$ is called {\bf $m$-sectorial} if it is sectorial and 
$\lambda I - A$ is surjective for some $\lambda \in \Ri$ with 
$\lambda < \gamma$.
Then an operator $A$ on $H$ is $m$-sectorial if and only if
$-A$ generates a holomorphic $C_0$-semigroup $S$, which is 
quasi-contractive on some sector, i.e.\ there exist 
$\theta \in (0,\frac{\pi}{2})$ and $\omega \in \Ri$ such that 
$\|e^{-\omega z} S_z\|_{\cl(H)} \leq 1$ for all $z \in \Sigma_\theta$.
(See \cite{Kat1} Theorem~IX.1.24 and \cite{Ouh5} proof of Theorem~1.58.)

The main theorem of this section is as follows.
We repeat that in our setting an operator is always univocal.

\begin{thm} \label{tsesqui201}
Let $H, V$ be Hilbert spaces and $j\colon V \to H$ a bounded linear operator 
such that $j(V)$ is dense in $H$.
Let $a \colon V \times V \to \Ci$ be a continuous sesquilinear form
which is $j$-elliptic.
Then one has the following.
\begin{tabel}
\item  \label{tsesqui201-1}
There exists an operator $A$ in $H$ such that 
for all $x,f \in H$ one has
$x \in D(A)$ and $Ax = f$ if and only if
\begin{quote}
there exists a $u \in V$ such that $j(u) = x$ and 
$a(u,v) = (f, j(v))_H$ for all $v \in V$.
\end{quote}
\item  \label{tsesqui201-2}
The operator $A$ of Statement~{\rm \ref{tsesqui201-1}} is $m$-sectorial.
\end{tabel}
\end{thm}
We call the operator $A$ in Statement~\ref{tsesqui201-1} of Theorem~\ref{tsesqui201}
the {\bf operator associated with} $(a,j)$.

In the proof of Theorem~\ref{tsesqui201} we need two subspaces of $V$
which we need throughout the paper.
Set
\[
D_H(a) = \{ u \in V : \mbox{there exists an }
 f \in H \mbox{ such that } a(u,v) = ( f, j(v))_H \mbox{ for all } v \in V \}
  \]
and 
\[
V(a) = \{ u \in V : a(u,v) = 0 \mbox{ for all } v \in \ker j \} 
.   \]
Clearly $D_H(a) \subset V(a)$ and $V(a)$ is closed in $V$.

\ruimte

\noindent
{\bf Proof of Theorem~\ref{tsesqui201}\hspace{5pt}}\
The proof consists of several steps.
\firststep\
First, we prove that the restriction map $j|_{V(a)} \colon V(a) \to H$ is injective.
If $u \in V(a)$ and $j(u) = 0$, then $a(u) = 0$.
The $j$-ellipticity (\ref{etsesqui201;1}) of $a$ then implies that $\|u\|_V = 0$.
So $u = 0$ and $j|_{V(a)}$ is injective.

\nextstep\
Next we prove Statement~\ref{tsesqui201-1}.
If $u \in V$, then it follows from the density of $j(V)$ in $H$
that there exists at most one $f \in H$ such that 
$a(u,v) = (f,j(v))_H$ for all $v \in V$.
But $j|{D_H(a)}$ is injective.
Therefore we can define the operator $A$ by $D(A) = j(D_H(a))$ and 
\begin{equation}
a(u,v) = (A j(u),j(v))_H
\;\;\;\;
\mbox{ for all } u \in D_H(a) \mbox{ and } v \in V .
\label{eSsesui2;1}
\end{equation}
(We emphasize that (\ref{eSsesui2;1}) is restricted to $u \in D_H(a)$
and need not to be valid for all $u \in V$ with $j(u) \in D(A)$.
An example will be given in Example~\ref{xsesqui323}.)

\nextstep\
Let $c$, $\omega$ and $\mu$ be as in (\ref{eSsesqui2.3}) and 
(\ref{etsesqui201;1}).
Let $x \in D(A)$. 
There exists a $u \in D_H(a)$ such that $x = j(u)$.
Then 
$((\omega I + A)x,x) = a(u) + \omega \, \|j(u)\|_H^2$ and 
$\RRe ((\omega I + A)x,x) \geq \mu \, \|u\|_V^2$.
Therefore
\begin{eqnarray*}
|\IIm ((\omega I + A)x,x)|
& = & |\IIm a(u)| 
\leq c \, \|u\|_V^2
\leq \frac{c}{\mu} \,  \RRe ((\omega I + A)x,x)
.  
\end{eqnarray*}
So $A$ is sectorial with vertex $-\omega$.

Finally, set $\lambda = \omega + 1$.
We shall show that the range of $\lambda I + A$ equals $H$.
Define the form $b$ on $V$ by 
$b(u,v) = a(u,v) + \lambda \, (j(u),j(v))_H$.
Then $b$ is continuous and coercive.
Let $f \in H$.
The Lax--Milgram theorem implies that there exists a unique $u \in V$ 
such that $b(u,v) = (f, j(v))_H$ for all $v \in V$.
Then $j(u) \in D(A)$ and $(\lambda I + A) j(u) = f$.
Thus $A$ is $m$-sectorial.
This proves Theorem~\ref{tsesqui201}.\hfill$\Box$

\ruimte

Although Theorem~\ref{tsesqui322} is a special case of 
Theorem~\ref{tsesqui340}, a short direct proof can be given 
at this stage.

\ruimte

\noindent
{\bf Proof of Theorem~\ref{tsesqui322}\hspace{5pt}}\
Denote by $V$ the completion of $(D(a),\| \cdot\|_a)$.
The injection of $(D(a),\| \cdot\|_a)$ into $H$ is continuous.
Hence there exists a $j\in \cl (V,H)$ such that $j(u)=u$ for all
$u \in D(a)$. 
Since $a$ is sectorial, there exists a unique 
continuous extension $\tilde a \colon V\times V \to \Ci$.
This extension is $j$-elliptic. 
Let $A$ be the operator associated with $(\tilde a, j)$.
If $u_1,u_2,\ldots \in D(a)$ with $u_1,u_2,\ldots$ convergent in $H$
 and $\RRe a(u_1), \RRe a(u_2),\ldots$ bounded, then 
$u_1,u_2,\ldots$ is bounded in $D(a)$. 
Therefore it has a weakly convergent subsequence in~$V$.
It follows from the density of $D(a)$ in $V$ that
$A$ equals the operator
from Theorem~\ref{tsesqui322}.
In particular, the operator is well-defined.
Now the result follows from Theorem~\ref{tsesqui201}.\hfill$\Box$

\ruimte

In the definition of $j$-elliptic the assumption is that (\ref{etsesqui201;1})
is valid for all $u \in V$. 
For a version of the Dirichlet-to-Neumann operator in Subsection~\ref{Ssesqui4.2}
this condition is too strong. 
One only needs (\ref{etsesqui201;1}) to be valid for all $u \in V(a)$
if in addition $V = V(a) + \ker j$ (cf.\ Theorem~\ref{tsesqui212}\ref{tsesqui212-1}).

\begin{cor} \label{csesqui201.2}
Let $H, V$ be Hilbert spaces and $j\colon V \to H$ a bounded linear operator 
such that $j(V)$ is dense in $H$.
Let $a \colon V \times V \to \Ci$ be a continuous sesquilinear form.
Suppose that there exist $\omega \in \Ri$
and $\mu > 0$
such that
\begin{equation}
\RRe a(u) + \omega \, \|j(u)\|_H^2 \geq \mu \, \|u\|_V^2
\label{ecsesqui201.2;1}
\end{equation}
for all $u\in V(a)$.
In addition suppose that $V = V(a) + \ker j$.
Then one has the following.
\begin{tabel}
\item  \label{csesqui201.2-1}
There exists an operator $A$ in $H$ such that 
for all $x,f \in H$ one has
$x \in D(A)$ and $Ax = f$ if and only if
\begin{quote}
there exists a $u \in V$ such that $j(u) = x$ and 
$a(u,v) = (f, j(v))_H$ for all $v \in V$.
\end{quote}
\item  \label{csesqui201.2-2}
The operator $A$ of Statement~{\rm \ref{csesqui201.2-1}} is $m$-sectorial.
\end{tabel}
\end{cor}
Again we call the operator $A$ in Statement~\ref{csesqui201.2-1} of Theorem~\ref{csesqui201.2}
the {\bf operator associated with} $(a,j)$.

\ruimte 

\noindent
{\bf Proof of Corollary~\ref{csesqui201.2}\hspace{5pt}}\
Define the form $b$ by $b = a|_{V(a) \times V(a)}$.
Then $(b,j|_{V(a)})$ satisfies the assumptions of Theorem~\ref{tsesqui201}.
Let $B$ be the operator associated with $(b,j|_{V(a)})$.

Clearly $D_H(a) \subset D_H(b)$.
Conversely, if $u \in D_H(b)$, then there is an $f \in H$ such that 
$b(u,v) = (f,j(v))_H$ for all $v \in V(a)$.
Then $a(u,v) = (f,j(v))_H$ for all $v \in V(a)$, but also 
for all $v \in \ker j$ by definition of $V(a)$.
Since $V = V(a) + \ker j$ by assumption, it follows that 
$u \in D_H(a)$.
Therefore $D_H(a) = D_H(b)$.
But $j|_{D_H(a)}$ is injective.
Hence the operator $B$ satisfies the requirements of 
Statement~\ref{csesqui201.2-1}.
Then Statement~\ref{csesqui201.2-2} is obvious.\hfill$\Box$

\ruimte

We return to the situation of Theorem~\ref{tsesqui201}.
If the form $a$ is $j$-elliptic and if $\tau \in \Ci$,
then obviously the operator $A + \tau I$ is associated with 
$(b,j)$, where $b$ is the $j$-elliptic form
$b(u,v) = a(u,v) + \tau \, (j(u),j(v))_H$ on $V$.

\ruimte

Although it is very convenient that we do not assume that the operator $j$ is 
injective, 
the second statement in the next proposition shows that without loss of generality
one might assume that 
$j$ is injective, by considering a different form.
The proposition is a kind of uniqueness result.
It determines the dependence of the operator on the choice of $V$.

\begin{prop} \label{psesqui204}
Suppose the form $a$ is $j$-elliptic and
let $A$ be the operator associated with $(a,j)$.
Then one has the following.
\begin{tabel}
\item \label{psesqui204-2}
If $U$ is a closed subspace of $V$ such that $D_H(a) \subset U$,
then $A$ equals the operator associated with $(a|_{U \times U}, j|_U)$.
If, in addition, the restriction $j|_U$ is injective, then $U = \overline{D_H(a)}$,
where the closure is taken in $V$.
\item \label{psesqui204-1}
$V(a) = \overline{D_H(a)}$.
Moreover, $j|_{V(a)}$ is injective
and $A$ equals the operator associated with $(a|_{V(a) \times V(a)},j|_{V(a)})$.
\item \label{psesqui204-3}
If $U$ is a closed subspace of $V(a)$ such that $j(U)$ is dense in $H$
and $A$ is the operator associated with $(a|_{U \times U}, j|_U)$, 
then $U = V(a)$.
\end{tabel}
\end{prop}
\proof\
`\ref{psesqui204-2}'.
Note that $j(U)$ and $j(V(a))$ both contain $j(D_H(a)) = D(A)$.
Therefore $j(U)$ and $j(V(a))$ are  
dense in $H$.
Let $b_1 = a|_{U \times U}$ and 
$b_2 = a|_{V(a) \times V(a)}$.
Further, let $B_1$ and $B_2$ be the operators associated with 
$(b_1,j|_U)$ and $(b_2,j|_{V(a)})$.
Then for all $u \in D_H(a)$ one deduces that 
$(A j(u), j(v))_H = a(u,v) = b_1(u,v)$ for all $v \in U$.
Therefore $u \in D_H(b_1)$ and $B_1 j(u) = A j(u)$.
So $A \subset B_1$.
But both $-A$ and $-B_1$ are semigroup generators.
Therefore $B_1 = A$.
Similarly, $A = B_2$.
Finally, if $j$ is injective on $U$,
then it follows from the inclusion $V(a) \subset U$ and 
the uniqueness theorem for closed sectorial forms,
\cite{Kat1} Theorem VI.2.7 that $U = V(a)$.
This proves Statement~\ref{psesqui204-2}.

`\ref{psesqui204-1}'.
The injectivity of $j|_{V(a)}$ has been proved in Step~1 of the 
proof of Theorem~\ref{tsesqui201}.
Then \ref{psesqui204-1} is a special case of \ref{psesqui204-2}.

Finally, Statement~\ref{psesqui204-3} follows from Statement~\ref{psesqui204-2}
with $a$ replaced by $a|_{U \times U}$.\hfill$\Box$

\ruimte

It is easy to construct examples with $V(a) \neq V$.
Therefore the injectivity condition in Proposition~\ref{psesqui204}\ref{psesqui204-2}
is necessary.

\begin{cor} \label{csesqui210}
Assume the notation and conditions of Corollary~{\rm \ref{csesqui201.2}}.
Let $A$ be the operator associated with $(a,j)$.
Then one has the following.
\begin{tabel}
\item \label{csesqui210-1}
$V(a) = \overline{D_H(a)}$.
Moreover, $j|_{V(a)}$ is injective
and $A$ equals the operator associated with $(a|_{V(a) \times V(a)},j|_{V(a)})$.
\item \label{csesqui210-3}
Let $U$ be a closed subspace of $V(a)$ such that $j(U)$ is dense in $H$.
Then $a|_{U \times U}$ is $j|_U$-elliptic.
Suppose $A$ is the operator associated with $(a|_{U \times U}, j|_U)$. 
Then $U = V(a)$.
\end{tabel}
\end{cor}
\proof\
Let $b = a|_{V(a) \times V(a)}$. 
Then it follows from the proof of Corollary~\ref{csesqui201.2} that 
$b$ is $j|_{V(a)}$-elliptic and $A$ is the associated operator.
Moreover, $D_H(a) = D_H(b)$.
Then 
\[
V(b)
= \{ u \in V(a) : a(u,v) = 0 \mbox{ for all } v \in V(a) \cap \ker j \} 
= V(a)
.   \]
It follows from Proposition~\ref{psesqui204}\ref{psesqui204-1} applied
to $(b,j|_{V(a)})$ that $V(b) = \overline{D_H(b)}$.
Hence $V(a) = \overline{D_H(a)}$.
Moreover, $(j|_{V(a)})|_{V(b)}$ is injective.
Therefore $j|_{V(a)}$ is injective.

If $U$ is a closed subspace of $V(a)$ such that $j(U)$ is dense in $H$,
then $a|_{U \times U}$ is $j|_U$-elliptic.
Then Statement~\ref{csesqui210-3} follows from 
Proposition~\ref{psesqui204}\ref{psesqui204-3}.\hfill$\Box$

\ruimte

In Subsection~\ref{Ssesqui4.2} we give an example that Proposition~\ref{psesqui204}\ref{psesqui204-2}
cannot be extended to the setting of Corollary~\ref{csesqui201.2}.

One can 
decompose a form $a$ in its real and imaginary parts as $a = h + i k$,
where $h,k \colon D(a) \times D(a) \to \Ci$ are symmetric sesquilinear forms.
We write $\Re a = h$ and $\Im a = k$.

The next theorem gives a connection between the current forms $a$ together with 
the map $j$
and the closed sectorial forms in Kato \cite{Kat1} Section~VI.2.

\begin{thm} \label{tsesqui212}
Suppose the form $a$ is $j$-elliptic and
let $A$ be the operator associated with $(a,j)$.
Then the following holds.
\begin{tabel}
\item \label{tsesqui212-1}
$\ker j \oplus V(a) = V$ as vector spaces.
\item \label{tsesqui212-2}
Let $a_c$ be the form on $H$ defined by
\[
D(a_c) = j(V) \mbox{ and } a_c(j(u),j(v)) 
= a(u,v) \quad (u,v \in V(a)).
\]
Then $a_c$ is the unique closed, sectorial form such that $A$ is 
associated with $a_c$.
\end{tabel}
\end{thm}
\proof\
`\ref{tsesqui212-1}'. 
Let $\omega \in \Ri$ and $\mu > 0$ be as in (\ref{etsesqui201;1}).
We can assume that $\omega = - 1$ and the form $a$ is coercive.  
Otherwise we replace $a$ by 
$(u,v) \mapsto a(u,v)+(\omega + 1) \, (j(u),j(v))_H$.
Let $h = \Re a$ and $k = \Im a$ be the real and imaginary part of $a$.
Then $\langle u,v\rangle := h(u,v)$ defines an equivalent scalar 
product on $V$.
So we may assume that $\|u\|_V = \|u\|_h$ for all $u \in V$.
Let $V_1 = \ker j$ and $V_2=(\ker j)^\perp$.
Moreover, let $\pi_1$ and $\pi_2$ be the projection from $V$
onto $V_1$ and $V_2$, respectively.
Then $h(u_1,v_2)=0$ for all
$u_1\in V_1$ and $v_2\in V_2$.
There exists a unique operator $T \in \cl(V)$
such that $k(u,v) = h(Tu,v)$ for all $u,v\in V$.
Let $T_{11} = \pi_1 \circ T|_{V_1} \in \cl(V_1)$ and 
$T_{12} = \pi_1 \circ T|_{V_2} \in \cl(V_2, V_1)$.
If $(u_1,u_2) \in V_1 \times V_2$, then
$u_1 + u_2 \in V(a)$ if and only if 
$0 = a(u_1 + u_2, v_1)
= h((I + i T)(u_1 + u_2),v_1) = h((I + i T_{11}) u_1 + i T_{12} u_2, v_1)$ for all
$v_1 \in V_1$.
So
\[
V(a)= \{ u_1 + u_2 : (u_1,u_2) \in V_1 \times V_2 
    \mbox{ and } (I + i T_{11}) u_1 + i T_{12} u_2 = 0 \} 
. \]
But $T$ is self-adjoint since $h$ and $k$ are symmetric.
So $I + i T_{11}$ is invertible.
Thus for all $u_2\in V_2$ there exists a $u_1\in V_1$ 
such that $u_1+u_2\in V(a)$.
Consequently, $j(V(a)) = j(V_2) = j(V)$. 
This implies that 
$\ker j + V(a) = V$. 
That the sum is direct has been proved in Step~1 of the proof of Theorem~\ref{tsesqui201}.

`\ref{tsesqui212-2}'. 
Define on $j(V(a))$ the scalar product carried over 
from $V(a)$ by $j$. 
Then the form~$a_c$ is clearly continuous
and elliptic, which is the same as sectorial and closed
(cf.\ Lemma~\ref{lsesqui241}).
The operator $A$ is clearly the operator associated with
$a_c$.\hfill$\Box$

\ruimte

We call the form $a_c$ in Theorem~\ref{tsesqui212}
the {\bf classical form} associated with $(a,j)$.
It equals the  classical form associated with the $m$-sectorial 
form $A$.
The proof of Theorem~\ref{tsesqui212} also allows to estimate the real part of the 
classical form of $a$ by the classical form of the real part of $a$ as follows.

\begin{prop} \label{psesqui214}
Suppose the form $a$ is $j$-elliptic and
let $A$ be the operator associated with $(a,j)$.
Suppose $\omega \leq -1$ in {\rm (\ref{etsesqui201;1})}.
Let $h$ be the real part of $a$ and 
$h_c$ the classical form associated with $(h,j)$.
Then $D(a_c) = D(h_c)$.
Moreover, there exists a constant $C > 0$ such that 
$\RRe a_c(x) \leq C \, h_c(x)$
for all $x \in j(V)$.
\end{prop}
\proof\
The first statement is obvious since $D(a_c) = j(V) = D(h_c)$.
We use the notation introduced in the proof of Theorem~\ref{tsesqui212}.
Moreover, we may assume that the inner product on $V$ is given by 
$(u,v) \mapsto h(u,v)$.
Let $u \in V(a)$. 
Then $(I + i T_{11}) u_1 + i T_{12} u_2 = 0$, where $u_1 = \pi_1(u)$ and $u_2 = \pi_2(u)$.
So $u_1 = - i (I + i T_{11})^{-1} T_{12} u_2$.
Moreover, $j(u) = j(u_2)$ and $u_2 \in V(h)$.
So $a_c(j(u)) = a(u)$ and $h_c(j(u)) = h_c(j(u_2)) = h(u_2) = \|u_2\|_V^2$.
Since the operator $(I + i T_{11})^{-1} T_{12}$ is bounded
one estimates
\begin{eqnarray*}
\RRe a_c(j(u))
& = & \RRe a(u)
= h(u)
= \|u_1\|_V^2 + \|u_2\|_V^2
\leq C \, \|u_2\|_V^2
= C \, h_c(j(u))
\end{eqnarray*}
where $C = \|(I + i T_{11})^{-1} T_{12}\|^2 + 1$.\hfill$\Box$

\ruimte

The next lemma gives a sufficient condition for the resolvents to be 
compact.

\begin{lemma} \label{lsesqui206}
Suppose the form $a$ is $j$-elliptic and
let $A$ be the operator associated with $(a,j)$.
If $j$ is compact, then $(\lambda I + A)^{-1}$ is compact for all 
$\lambda \in \Ci$ with $\RRe \lambda > \omega$, where $\omega$ is as
in {\rm (\ref{etsesqui201;1})}.
\end{lemma}
\proof\
By the Lax--Milgram theorem there exists a $B \in \cl(H,V)$ such that 
\[
(f,j(v))_H = a(Bf, v) + \lambda \, (j(Bf), j(v))_H
\]
for all $f \in H$ and $v \in V$.
Then $B(H) \subset D_H(a)$ and 
$(A + \lambda I) j(Bf) = f$ for all $f \in H$.
Therefore $(\lambda I + A)^{-1} = j \circ B$ is compact.\hfill$\Box$

\begin{remarkn} \label{rsesqui245}
If $B$ is the operator associated with $(a^*,j)$ where
$a^*$ is the $j$-elliptic form on $V$ given by $a^*(u,v) = \overline{a(v,u)}$,
then $A^*$ is an extension of $B$. 
But both $-A^*$ and $-B$ are generators of semigroups.
Therefore $A^*$ is the operator associated with $(a^*,j)$.
\end{remarkn}

In \cite{Ouh5} Theorem~2.2 there is a characterization of 
closed convex subsets which are invariant under the semigroup $S$.
For a background of this theorem we refer to the Notes for Section~2.1 in \cite{Ouh5}.
Using the two statements of Theorem~\ref{tsesqui212}, the theorem
of Ouhabaz can be reformulated in the current context.
Recall that a sesquilinear form $b$ is called {\bf accretive}
if $\RRe b(u) \geq 0$ for all $u \in D(b)$.

\begin{prop} \label{psesqui202}
Suppose the form $a$ is $j$-elliptic,
let $A$ be the operator associated with $(a,j)$ and $S$ the 
semigroup generated by $-A$.
Moreover, suppose that $a$ is accretive.
Let $C \subset H$ be a non-empty closed convex set and $P \colon H \to C$ the 
orthogonal projection.
Then the following conditions are equivalent.
\begin{tabeleq}
\item \label{psesqui202-1}
$S_t C \subset C$ for all $t > 0$.
\item \label{psesqui202-2}
For all $u \in V$ there exists a $w \in V$ such that 
\[
Pj(u) = j(w)
\;\;\;\; \mbox{and} \;\;\;\;
\RRe a(w, u-w) \geq 0 
.  \]
\item \label{psesqui202-3}
For all $u \in V$ there exists a $w \in V$ such that 
\[
Pj(u) = j(w)
\;\;\;\; \mbox{and} \;\;\;\;
\RRe a(u, u-w) \geq 0 
.  \]
\item \label{psesqui202-4}
There exists a dense subset $D$ of $V$ such that 
for all $u \in D$ there exists a $w \in V$ such that 
\[
Pj(u) = j(w)
\;\;\;\; \mbox{and} \;\;\;\;
\RRe a(w, u-w) \geq 0 
.  \]
\end{tabeleq}
\end{prop} 
\proof\
`\ref{psesqui202-1}$\Rightarrow$\ref{psesqui202-2}'.
Let $u \in V$.
By Theorem~\ref{tsesqui212} there exists a $u' \in V(a)$ such that 
$j(u') = j(u)$.
Then $P j(u') \in D(a_c)$ by \cite{Ouh5} Theorem~2.2 1)$\Rightarrow$2).
So there exists a $w \in V(a)$ such that $P j(u') = j(w)$.
Then 
$\RRe a(w,u'-w)
= \RRe a_c(j(w), j(u') - j(w))
= \RRe a_c(P j(u'), j(u') - P j(u'))
\geq 0$
again by \cite{Ouh5} Theorem~2.2 1)$\Rightarrow$2).
But $a(w,u-u') = 0$ since $w \in V(a)$ and $u-u' \in \ker j$.
So $\RRe a(w, u-w) \geq 0$.

`\ref{psesqui202-2}$\Rightarrow$\ref{psesqui202-3}'.
Trivial, since $\RRe a(u-w, u-w) \geq 0$.

`\ref{psesqui202-3}$\Rightarrow$\ref{psesqui202-1}'.
Let $u \in V(a)$.
By assumption there exists a $w \in V$ such that $P j(u) = j(w)$ and 
$\RRe a(u, u-w) \geq 0$.
Let $w' \in V(a)$ be such that $j(w) = j(w')$.
Then $a(u,w-w') = 0$ since $u \in V(a)$ and $w-w' \in \ker j$.
So $\RRe a(u, u-w') \geq 0$ and $\RRe a_c(j(u), j(u) - P j(u)) \geq 0$.
Then the implication follows from \cite{Ouh5} Theorem~2.2 3)$\Rightarrow$1).

`\ref{psesqui202-2}$\Rightarrow$\ref{psesqui202-4}'.
Trivial.

`\ref{psesqui202-4}$\Rightarrow$\ref{psesqui202-2}'.
Since $a$ is continuous there exists a $c > 0$ such that 
$|a(u,v)| \leq c \, \|u\|_V \, \|v\|_V$ for all $u,v \in V$.
Let $u \in V$.
There exist $u_1,u_2,\ldots \in D$ such that $\lim u_n = u$ in $V$.
For all $n \in \Ni$ there exists by assumption a $w_n \in V$ such that 
$Pj(u_n) = j(w_n)$ and $\RRe a(w_n, u_n - w_n) \geq 0$.
Let $\mu$ and $\omega$ be as in (\ref{etsesqui201;1}). 
Then 
\begin{eqnarray*}
\mu \, \|w_n\|_V^2
& \leq & \RRe a(w_n) + \omega \, \|j(w_n)\|_H^2  \\
& = & \RRe a(w_n,u_n) - \RRe a(w_n, u_n - w_n) + \omega \, \|j(w_n)\|_H^2 \\
& \leq & \RRe a(w_n,u_n) + \omega \, \|j(w_n)\|_H^2 \\ 
& \leq & c \, \|w_n\|_V \, \|u_n\|_V + \omega \, \|P j(u_n)\|_H^2  
\end{eqnarray*}
for all $n \in \Ni$.
Since $ \{ u_n : n \in \Ni \} $ is bounded in $V$ 
and $ \{ Pj(u_n) : n \in \Ni \} $ is bounded in $H$ by continuity of 
$j$ and contractivity of $P$, it follows that the set 
$ \{ w_n : n \in \Ni \} $ is bounded in $V$.
So there exist $w \in V$ and a subsequence $w_{n_1},w_{n_2},\ldots$
of $w_1,w_2,\ldots$ such that $\lim_{k \to \infty} w_{n_k} = w$ weakly in $V$.
Then $\lim_{k \to \infty} Pj(u_{n_k}) = \lim j(w_{n_k}) = j(w)$
weakly in $H$.
On the other hand, the continuity of $j$ and $P$ gives 
$\lim_{n \to \infty} Pj(u_n) = Pj(u)$ strongly in $H$.
So $Pj(u) = j(w)$.
Since $\RRe a(w_n, u_n - w_n) \geq 0$ one has
$\RRe a(w_n) \leq \RRe a(w_n,u_n)$ for all $n \in \Ni$.
Moreover, 
$\lim_{k \to \infty} \RRe a(w_{n_k},u_{n_k}) = \RRe a(w,u)$.
In addition, since $a$ is accretive and $j$-elliptic
it follows that $v \mapsto ( \RRe a(v) + \varepsilon \, \|j(v)\|_H^2)^{1/2}$
is an equivalent norm associated with an inner product on $V$ for all $\varepsilon > 0$.
Therefore $\RRe a(w) \leq \liminf_{k \to \infty} \RRe a(w_{n_k})$.
So $\RRe a(w) \leq \RRe a(w,u)$ and $\RRe a(w, u-w) \ge 0$ as required.\hfill$\Box$

\section{Generation theorems in the incomplete case} \label{Ssesqui3}

In this section we consider forms for which the form domain is 
not necessarily a Hilbert space.
First we reformulate the complete case.

Let $a \colon D(a) \times D(a) \to \Ci$ be a sesquilinear form, 
where the domain $D(a)$ of $a$ is a complex vector
space, the domain of~$a$.
Let $H$ be a Hilbert space and $j \colon D(a) \to H$ a linear map.
We say that $a$ is a {\bf $j$-sectorial form} if 
there are $\gamma \in \Ri$, called a {\bf vertex},
and $\theta \in [0,\frac{\pi}{2})$, called a {\bf semi-angle}, such that 
\[
a(u) - \gamma \, \|j(u)\|_H^2
\in \Sigma_\theta
\]
for all $u \in D(a)$.
If $a$ is $j$-sectorial with vertex $\gamma$, then we 
define the semi-inner product $(\, \cdot , \cdot \, )_a$ in the 
space $D(a)$ by 
\[
(u, v)_a
= (\Re a)(u,v) + (1-\gamma) \, (j(u) , j(v))_H
.  \]
Again we do not include the $\gamma$ in the notation.
Then the associated seminorm $\|\cdot\|_a$ is a norm if and only if 
$\RRe a(u) = j(u) = 0$ implies
$u = 0$ for all $u \in D(a)$.
A $j$-sectorial form $a$ is called {\bf closed} if $\|\cdot\|_a$ is a norm and 
$(D(a), \|\cdot\|_a)$ is a Hilbert space.
This term coincides with the term for closed forms in \cite{Kat1}
Section~VI.1.3 if $j$ is an inclusion map.

The alluded reformulation is as follows.

\begin{lemma} \label{lsesqui241}
Let $V$ be a vector space,
$a \colon V \times V \to \Ci$ a sesquilinear form,
$H$ a Hilbert space and $j \colon V \to H$ a linear map.
Then the following are equivalent.
\begin{tabeleq}
\item \label{lsesqui241-1}
The form $a$ is $j$-sectorial and closed.
\item \label{lsesqui241-2}
There exists a norm $\|\cdot\|_V$ on $V$ such that $V$ is a Banach  space,
the map $j$ is bounded from $(V,\|\cdot\|_V)$ into $H$, 
the form $a$ is $j$-elliptic and $a$ is continuous.
\end{tabeleq}
Moreover, if Condition~{\rm \ref{lsesqui241-2}} is valid, then the norms $\|\cdot\|_a$ and
$\|\cdot\|_V$ are equivalent.
\end{lemma}
\proof\
The easy proof is left to the reader.\hfill$\Box$

\ruimte

In this section we drop the assumption that $(D(a), \|\cdot\|_a)$ is complete. 
So $H$ is a Hilbert space, $a \colon D(a) \times D(a) \to \Ci$ is a sesquilinear form,
$j \colon D(a) \to H$ is a linear map and we assume that 
$a$ is merely $j$-sectorial and $j(D(a))$ is dense in $H$.
We will again associate a sectorially bounded holomorphic semigroup
generator with~$(a,j)$.
The next theorem is an extension of both Theorem~\ref{tsesqui322} and 
Theorem~\ref{tsesqui201}.

In a natural way one can define the notion of Cauchy sequence in a 
semi-normed vector space.

\begin{thm} \label{tsesqui340}
Let $a$ be a sesquilinear form, 
$H$ a Hilbert space and $j \colon D(a) \to H$ a linear map.
Assume that $a$ is $j$-sectorial and $j(D(a))$ is dense in $H$.
Then one has the following.
\begin{tabel}
\item  \label{tsesqui340-1}
There exists an operator $A$ in $H$ such that 
for all $x,f \in H$ one has
$x \in D(A)$ and $Ax = f$ if and only if there exist $u_1,u_2,\ldots \in D(a)$
such that
\begin{tabelR}
\item \label{tsesqui340-11}
$\lim_{n \to \infty} j(u_n) = x$ weakly in $H$,
\item \label{tsesqui340-12}
$\sup_{n \in \Ni} \RRe a(u_n) < \infty$, and, 
\item \label{tsesqui340-13}
$\lim_{n \to \infty} a(u_n,v) = (f,j(v))_H$ for all $v \in D(a)$. 
\end{tabelR}
\item  \label{tsesqui340-2}
The operator $A$ of Statement~{\rm \ref{tsesqui340-1}} is $m$-sectorial.
\item  \label{tsesqui340-3}
Let $x,f \in H$.
Then $x \in D(A)$ and $Ax = f$
if and only if 
there exists a Cauchy sequence $u_1,u_2,\ldots$ in $D(a)$
such that $\lim_{n \to \infty} j(u_n) = x$ in $H$ and 
$\lim_{n \to \infty} a(u_n,v) = (f,j(v))_H$ for all $v \in D(a)$.
\end{tabel}
\end{thm}

If $V_0$ is a vector space with a semi-inner product, then there exist a 
Hilbert space $V$ and an isometric map $q \colon V_0 \to V$ 
such that $q(V_0)$ is dense in $V$.
Then $V$ and $q$ are unique, up to unitary equivalence.
We call $(V,q)$ the {\bf completion} of $V_0$.
If also $V_0'$ is a vector space with a semi-inner product
and $(V',q')$ is its completion, then for every linear map
$T_0 \colon V_0 \to V_0'$ and $c \geq 0$ such that 
$\|T_0 u\|_{V_0'} \leq c \, \|u\|_{V_0}$ for all $u \in V_0$,
there exists a unique $T \in \cl(V,V')$ such that 
$q' \circ T_0 = T \circ q$. 
Then $\|T\| \leq c$.
We call $T$ the {\bf continuous extension} of $T_0$ to~$V$.

\ruimte

\noindent
{\bf Proof of Theorem~\ref{tsesqui340}\hspace{5pt}\ }\
Let $(V,q)$ be the completion of $D(a)$.
Since $\|j(u)\|_H \leq \|u\|_a$ for all $u \in D(a)$
the continuous extension $\tilde j \in \cl(V,H)$ of $j$ is a 
contraction.
Note that $\tilde j \circ q = j$ and $\tilde j(V)$ is dense in $H$.
Next, since 
\[
|a(u,v) - \gamma \, (j(u), j(v))_H|
\leq (1 + \tan \theta) \, \|u\|_a \, \|v\|_a
\]
for all $u,v \in D(a)$, where $\theta$ is the semi-angle of $a$
and we used the estimate (1.15) of Subsection~VI.1.2 in \cite{Kat1},
there exists a unique continuous sesquilinear form 
$\tilde a \colon V \times V \to \Ci$
such that $\tilde a(q(u), q(v)) = a(u,v)$ for all $u,v \in D(a)$.
Then $\tilde a$ is $\tilde j$-sectorial
with vertex $\gamma$ and semi-angle $\theta$.
Moreover, $\tilde a$ is $\tilde j$-elliptic.
Now let $A$ be the 
operator associated with $(\tilde a, \tilde j)$.

Let $x,f \in H$.
We next show that the statements
\begin{tabeleq}
\item  \label{tsesqui340=1}
$x \in D(A)$ and $Ax = f$,
\item  \label{tsesqui340=2}
there exists a Cauchy sequence $u_1,u_2,\ldots$ in $(D(a), \|\cdot\|_a)$
such that $\lim j(u_n) = x$ and 
$\lim_{n \to \infty} a(u_n,v) = (f,j(v))_H$ for all $v \in D(a)$, and
\item  \label{tsesqui340=3}
there exists a bounded sequence $u_1,u_2,\ldots$ in $(D(a), \|\cdot\|_a)$
such that $\lim j(u_n) = x$ weakly in $H$ and 
$\lim_{n \to \infty} a(u_n,v) = (f,j(v))_H$ for all $v \in D(a)$
\end{tabeleq}
are equivalent.

`\ref{tsesqui340=1}$\Rightarrow$\ref{tsesqui340=2}'.
It follows from the definition (\ref{eSsesui2;1}) that there exists 
a $\tilde u \in V$ such that $\tilde j(\tilde u) = x$ and
$\tilde a(\tilde u,\tilde v) = (f,\tilde j(\tilde v))_H$
for all $\tilde v \in V$.
Then there exist $u_1,u_2,\ldots \in D(a)$ such that 
$\lim q(u_n) = \tilde u$ in $V$.
Hence $u_1,u_2,\ldots$ is a Cauchy sequence in $(D(a), \|\cdot\|_a)$.
Moreover,  
\[
(f,j(v))_H
= (f,\tilde j(q(v)))_H
= \tilde a(\tilde u, q(v))
= \lim \tilde a(q(u_n), q(v))
= \lim a(u_n,v)
\]
for all $v \in D(a)$ and 
$\lim j(u_n) = \lim \tilde j(q(u_n)) = \tilde j(\tilde u) = x$ in~$H$.

`\ref{tsesqui340=2}$\Rightarrow$\ref{tsesqui340=3}'. 
Trivial.

`\ref{tsesqui340=3}$\Rightarrow$\ref{tsesqui340=1}'.
Since $q(u_1),q(u_2),\ldots$ is a bounded sequence in $V_0$ the
weak limit $\tilde u = \lim q(u_n)$ exists in~$V$ after passing to a 
subsequence, if necessary.
Then $\tilde j(\tilde u) = \lim \tilde j(q(u_n)) = \lim j(u_n) = x$ weakly in~$H$.
Moreover, 
\[
\tilde a(\tilde u,q(v))
= \lim \tilde a(q(u_n), q(v))
= \lim a(u_n,v)
= (f,j(v))_H
= (f,\tilde j(q(v)))_H
\]
for all $v \in D(a)$.
Since $q(D(a))$ is dense in $V$ one deduces that 
$\tilde a(\tilde u,\tilde v) = (f,\tilde j(\tilde v))_H$
for all $\tilde v \in V$.
So $x \in D(A)$ and $A x = f$ as required.

We have proved the existence of the operator $A$ in Statement~\ref{tsesqui340-1}
of the theorem, together with the characterization~\ref{tsesqui340-3}.
Now Statement~\ref{tsesqui340-2} follows from 
Theorem~\ref{tsesqui201}.\hfill$\Box$

\ruimte

We call the operator $A$ in Statement~\ref{tsesqui340-1} of 
Theorem~\ref{tsesqui340} the {\bf operator associated with} $(a,j)$.
Note that this is the same operator as in Theorem~\ref{tsesqui201}
if $D(a)$ was provided with a Hilbert space
structure such that $j$ is continuous, $a$ is continuous and $a$ is $j$-elliptic.

\ruimte

In the proof of Theorem~\ref{tsesqui340} we also proved the following fact.

\begin{prop} \label{psesqui340.5}
Let $a$ be a sesquilinear form, 
$H$ a Hilbert space and $j \colon D(a) \to H$ a linear map.
Assume that $a$ is $j$-sectorial and $j(D(a))$ is dense in $H$.
Let $(V,q)$ be the completion of $D(a)$.
Then there exists a unique continuous sesquilinear form 
$\tilde a \colon V \times V \to \Ci$ such that 
$\tilde a(q(u),q(v)) = a(u,v)$ for all $u,v \in D(a)$.
Moreover, $\tilde a$ is $\tilde j$-elliptic, where 
$\tilde j$ is the continuous extension of $j$ to $V$
and the operator associated with $(a,j)$ equals the operator 
associated with $(\tilde a,\tilde j)$.
\end{prop}

\begin{remarkn} \label{rsesqui345}
Let $a$ be a sesquilinear form, 
$H$ a Hilbert space and $j \colon D(a) \to H$ a linear map.
Suppose that $a$ is $j$-sectorial.
Let $D$ be {\bf core} for $a$, i.e.\ a dense subspace of $D(a)$.
Then $j(D)$ is dense in $H$ and the operator associated with $(a,j)$ equals the 
operator associated with $(a|_{D \times D}, j|_D)$.
This follows immediately from 
Theorem~\ref{tsesqui340}\ref{tsesqui340-3}.
\end{remarkn}

\begin{remarkn} \label{rsesqui346}
Let $a$ be a sesquilinear form, 
$H$ a Hilbert space and $j \colon D(a) \to H$ a linear map.
Assume that $a$ is $j$-sectorial and $j(D(a))$ is dense in $H$.
Then $a^*$ is $j$-sectorial.
Moreover, if $B$ is the operator associated with $(a^*,j)$
and $A$ is the operator associated with $(a,j)$,
then $B = A^*$.
Indeed, using the notation as in the proof of Theorem~\ref{tsesqui340}
it follows that $A$ is the operator associated with 
$(\tilde a, \tilde j)$.
Moreover, $(u,v)_{a^*} = (u,v)_a$
for all $u,v \in D(a) = D(a^*)$.
Therefore $(V,q)$ is also the completion of $D(a^*)$.
Then $\widetilde{a^*} = (\tilde a)^*$.
By construction the operator $B$ is the operator associated with 
$(\widetilde{a^*},\tilde j) = ((\tilde a)^*,\tilde j)$.
Hence $B = A^*$ by Remark~\ref{rsesqui245}.
In particular, if $a$ is symmetric, then $A$ is self-adjoint.
\end{remarkn}

\begin{remarkn} \label{rsesqui346.2}
It follows from the construction that the operator $\lambda I + A$ is invertible for all 
$\lambda > (-\gamma) \vee 0$ if $A$ is the operator associated with a 
$j$-sectorial form $a$ with vertex $\gamma$. 
\end{remarkn}

The next theorem is of the nature of
\cite{Kat1} Theorem~VIII.3.6.
If $F_1,F_2,\ldots$ are subsets of a set $F$, then define
$\liminf_{n \to \infty} F_n = \bigcup_{n=1}^\infty \bigcap_{k=n}^\infty F_k$.

\begin{thm} \label{tsesqui499}
Let $a$ be a sesquilinear form, 
$H$ a Hilbert space and $j \colon D(a) \to H$ a linear map.
Assume that $a$ is $j$-sectorial with vertex $\gamma$.
For all $n \in \Ni$ let $a_n$ be a sesquilinear form with 
$D(a_n) \subset D(a)$.
Suppose that there exist $\theta \in [0,\frac{\pi}{2})$
and for all $n \in \Ni$ a $\gamma_n \in \Ri$ such that 
\begin{equation}
a_n(u) - a(u) - \gamma_n \, \|j(u)\|_H^2 
\in \Sigma_\theta
\label{etsesqui499;1}
\end{equation}
for all $u \in D(a_n)$.
Assume that $\lim_{n \to \infty} \gamma_n = 0$.
Moreover, suppose that there exists a core $D$ for $a$ 
such that $D \subset \liminf_{n \to \infty} D(a_n)$ and 
$\lim_{n \to \infty} a_n(u) = a(u)$ for all 
$u \in D$.
Finally, suppose that $j(D(a_n))$ is dense in $H$ for all $n \in \Ni$.
Let $A$ be the operator associated with $(a,j)$ and for all 
$n \in \Ni$ let $A_n$ be the operator associated with $(a_n,j|_{D(a_n)})$.
Fix $\lambda > (-\gamma) \vee 0$.
Then 
\[
\lim_{n \to \infty} (\lambda I + A_n)^{-1} f = (\lambda I + A)^{-1} f
\]
for all $f \in H$.
\end{thm}
\proof\
Without loss of generality we may assume that $\gamma = 0$
and $\gamma_n < 1$ for all $n \in \Ni$.
Then $a_n$ is $j$-sectorial with vertex $\gamma_n$
and $D(a_n)$ has the norm
$\|u\|_{a_n}^2 = \RRe a_n(u) + (1-\gamma_n) \, \|j(u)\|_H^2$.
We use the construction as in the proof of Theorem~\ref{tsesqui340}.
For the form $a$ we construct $V$, $q$, $\tilde j$, $\tilde a$
and for the form $a_n$ we construct $V_n$, $q_n$, $\tilde j_n$,
$\tilde a_n$.

Let $n \in \Ni$.
It follows from (\ref{etsesqui499;1}) that 
$\|u\|_a^2 \leq \|u\|_{a_n}^2$ for all $u \in D(a_n)$.
Let $\Phi_n$ be the continuous extension of the inclusion map
$D(a_n) \subset D(a)$.
So $\Phi_n \in \cl(V_n,V)$ and 
$\Phi_n \circ q_n = q$.
Then
$\tilde j_n(q_n(u)) = j(u) = \tilde j(q(u)) = \tilde j(\Phi_n(q_n(u)))$
for all $u \in D(a_n)$ and by density
$\tilde j_n = \tilde j \circ \Phi_n$.
Define the sectorial form $b_n \colon D(a_n) \times D(a_n) \to \Ci$
by 
\[
b_n(u,v) = a_n(u,v) - a(u,v) - \gamma_n \, (j(u),j(v))_H
.  \]
Then $|b_n(u)| \leq \|u\|_{a_n}^2$, so there
exists a unique continuous accretive sectorial form 
$\tilde b_n \colon V_n \times V_n \to \Ci$ such that 
$\tilde b_n(q_n(u),q_n(v)) = b_n(u,v)$ for all $u,v \in D(a_n)$.
Then 
\begin{equation}
\tilde a_n(u,v) 
= \tilde a(\Phi_n(u),\Phi_n(v)) 
   + \tilde b_n(u,v)
   + \gamma_n \, (\tilde j_n(u),\tilde j_n(v))_H
\label{etsesqui499;2}
\end{equation}
first for all $u,v \in q_n(D(a_n))$ and then by density for all 
$u,v \in V_n$.

In order not to duplicate too much of the proof for the current theorem
for the proof of Theorem~\ref{tsesqui402}
we first prove a little bit more.
Let $f,f_1,f_2,\ldots \in H$ and suppose that $\lim f_n = f$ 
weakly in $H$.
For all $n \in \Ni$ there exists a unique $\tilde u_n \in D_H(\tilde a_n)$ such that 
$\tilde j_n(\tilde u_n) = (\lambda I + A_n)^{-1} f_n$.
Set $u_n = \Phi_n(\tilde u_n) \in V$.
Then $\tilde j(u_n) = \tilde j_n(\tilde u_n)$.
We shall show that there exists a subsequence $(u_{n_k})$ of $(u_n)$
and a $u \in D_H(\tilde a)$
such that $\lim u_{n_k} = u$ weakly in $V$ and 
$\tilde j(u) = (\lambda I + A)^{-1} f$.

Since $\tilde u_n \in D_H(\tilde a_n)$ and 
$\lambda \, \tilde j_n(\tilde u_n) + A_n \tilde j_n(\tilde u_n) = f_n$
it follows from (\ref{eSsesui2;1}) that 
\begin{equation}
\lambda \, (\tilde j_n(\tilde u_n), \tilde j_n(v))_H
   + \tilde a_n(\tilde u_n,v)
= (f_n, \tilde j_n(v))_H
\label{etsesqui499;3}
\end{equation}
for all $v \in V_n$.
Taking $v = \tilde u_n$ in (\ref{etsesqui499;3})  and using (\ref{etsesqui499;2})
we obtain
\begin{eqnarray}
(\lambda + \gamma_n) \|\tilde j_n(\tilde u_n)\|_H^2 + \RRe \tilde a(u_n) 
     +  \RRe \tilde b_n(\tilde u_n) 
& = & \RRe (f_n,\tilde j_n(\tilde u_n))_H  \label{etsesqui499;4}  \\
& \leq & \|f_n\|_H \, \|\tilde j_n(\tilde u_n)\|_H 
\leq \tfrac{\lambda}{2} \, \|\tilde j_n(\tilde u_n)\|_H^2 
         + \tfrac{2}{\lambda} \|f_n\|_H^2
   \nonumber
.
\end{eqnarray}
Since $\frac{\lambda}{4} + \gamma_n \geq 0$ for large $n$
this implies that the set $ \{ \tilde j(u_n) : n \in \Ni \} 
= \{ \tilde j_n(\tilde u_n) : n \in \Ni \} $ is bounded in $H$,
and that the two sets 
$ \{ \RRe \tilde a(u_n) : n \in \Ni \} $
and $ \{ \RRe \tilde b_n(\tilde u_n) : n \in \Ni \} $
are bounded.
In particular the sequence $u_1,u_2,\ldots$ is bounded in $V$.
Passing to a subsequence, if necessary, it follows that there exists a 
$u \in V$ such that $\lim u_n = u$ weakly in $V$.
Then $\lim \tilde j(u_n) = \tilde j(u)$ weakly in $H$.

Let $n \in \Ni$.
Then $\tilde b_n$ is $\tilde j_n$-sectorial with vertex $0$ and semi-angle~$\theta$.
Therefore 
\[
|\tilde b_n(\tilde u_n,v)|
\leq (1 + \tan \theta) \Big( \RRe \tilde b_n(\tilde u_n) \Big)^{1/2} 
          \Big( \RRe \tilde b_n(v) \Big)^{1/2} 
\]
for all $v \in V_n$.
Now let $v \in D$.
Then $\lim_{n \to \infty} \RRe b_n(v) = 0$ by assumption.
Hence
$\lim_{n \to \infty} \tilde b_n(\tilde u_n,q_n(v)) = 0$.
It follows from (\ref{etsesqui499;2}) and (\ref{etsesqui499;3}) that 
\[
\lambda \, (\tilde j(u_n),j(v))_H
   + \tilde a(u_n,q(v))
   + \tilde b_n(\tilde u_n,q_n(v)) 
   + \gamma_n \, (\tilde j(u_n),j(v))_H
= (f_n,j(v))_H
.  \]
Taking the limit $n \to \infty$ gives
\begin{equation}
\lambda \, (\tilde j(u),j(v))_H + \tilde a(u,q(v))
= (f,j(v))_H
\label{etsesqui499;5}
\end{equation}
for all $v \in D$.
Since $D$ is a core for $a$ one deduces that (\ref{etsesqui499;5}) is valid for 
all $v \in D(a)$ and then again by density one establishes that 
\begin{equation}
\lambda \, (\tilde j(u),\tilde j(v))_H + \tilde a(u,v)
= (f,\tilde j(v))_H
\label{etsesqui499;6}
\end{equation}
for all $v \in V$.
Thus $u \in D_H(\tilde a)$, and by definition of $A$ it follows that 
$\tilde j(u) = (\lambda I + A)^{-1} f$.

Now we prove the theorem.
Let $f \in H$ and apply the above with $f_n = f$ for all $n \in \Ni$.
In order to deduce that $\lim \tilde j(u_n) = \tilde j(u)$ strongly in~$H$, 
by Proposition~3.6 in \cite{HirL}
it suffices to show that $\limsup \|\tilde j(u_n)\|_H \leq \|\tilde j(u)\|_H$.

Substituting $v = u_n$ in (\ref{etsesqui499;6}) gives 
\[
\lambda \, (\tilde j(u), \tilde j(u_n))_H + \tilde a(u,u_n) = (f,\tilde j(u_n))_H
\]
for all $n \in \Ni$.
Hence by (\ref{etsesqui499;4}) one deduces that 
\begin{eqnarray*}
\lambda \|\tilde j(u_n)\|_H^2 
& \leq & \lambda \, \|\tilde j_n(\tilde u_n)\|_H^2 
     + \RRe \tilde b_n(\tilde u_n)   \\
& = & \RRe \Big( (f,\tilde j(\tilde u_n))_H - \tilde a(u_n) \Big) 
    - \gamma_n \, \|\tilde j(u_n)\|_H^2  \\
& = & \RRe \Big( \lambda \, (\tilde j(u), \tilde j(u_n))_H + 
                    \tilde a(u,u_n) - \tilde a(u_n) \Big) 
    - \gamma_n \, \|\tilde j(u_n)\|_H^2
\end{eqnarray*}
for all $n \in \Ni$.
But $\RRe \tilde a(u) \leq \liminf \RRe \tilde a(u_n)$ by \cite{Kat1}, Lemma~VIII.3.14a.
Therefore $\limsup \lambda \, \|\tilde j(u_n)\|_H^2 \leq \RRe \lambda \, \|\tilde j(u)\|_H^2
= \lambda \, \|\tilde j(u)\|_H^2$ and the strong convergence follows.

We have shown that there exists a subsequence 
$n_1,n_2,\ldots$ of the sequence $1,2,\ldots$ such that 
$\lim_{k \to \infty} (\lambda I + A_{n_k})^{-1} f = (\lambda I + A)^{-1} f$.
But this implies that 
\[
\lim_{n \to \infty} (\lambda I + A_n)^{-1} f = (\lambda I + A)^{-1} f
\]
and the proof of the theorem is complete.\hfill$\Box$

\ruimte

For compact maps one obtains a stronger convergence in Theorem~\ref{tsesqui499}.

\begin{thm} \label{tsesqui402}
Assume the notation and conditions of Theorem~{\rm \ref{tsesqui499}}.
Suppose in addition that the map $j \colon D(a) \to H$ is compact.
Then 
\[
\lim_{n \to \infty} \|(\lambda I + A_n)^{-1} - (\lambda I + A)^{-1}\| = 0
\]
for all $\lambda > (-\gamma) \vee 0$.
\end{thm}
\proof\
Suppose not. 
Then there exist $\varepsilon > 0$, $n_1,n_2,\ldots \in \Ni$ 
and $f_1,f_2,\ldots \in H$ such that $n_k < n_{k+1}$, $\|f_k\|_H \leq 1$
and $\|(\lambda I + A_{n_k})^{-1} f_k - (\lambda I + A)^{-1} f_k\| \geq \varepsilon$
for all $k \in \Ni$.
Passing to a subsequence, if necessary, there exists an $f \in H$ 
such that $\lim_{k \to \infty} f_{n_k} = f$ weakly in $H$.
For all $k \in \Ni$ there exists a $\tilde u_k \in D_H(\tilde a_{n_k})$ such that 
$\tilde j_{n_k}(\tilde u_k) = (\lambda I + A_{n_k})^{-1} f_k$.
Let $u_k = \Phi_{n_k}(\tilde u_k)$, where we use the notation
as in the proof of Theorem~\ref{tsesqui499}.
Then it follows from the first part of the proof of Theorem~\ref{tsesqui499}
that there exists a $u \in D_H(\tilde a)$ such that, 
after passing to a subsequence if necessary, 
$\lim_{k \to \infty} u_k = u$ weakly in $V$ and 
$\tilde j(u) = (\lambda I + A)^{-1} f$.
Since $j$ is compact, the map $\tilde j$ is compact.
Therefore
\[
\lim_{k \to \infty}  (\lambda I + A_{n_k})^{-1} f_k
= \lim \tilde j(u_k)
= \tilde j(u)
= (\lambda I + A)^{-1} f
\]
strongly in $H$.
Moreover, 
$\lim_{k \to \infty} (\lambda I + A)^{-1} f_k = (\lambda I + A)^{-1} f$
by Lemma~\ref{lsesqui206}.
So $\lim_{k \to \infty} \|(\lambda I + A_{n_k})^{-1} f_k - (\lambda I + A)^{-1} f_k\| = 0$.
This is a contradiction.\hfill$\Box$

\ruimte

If $a$ is symmetric and $j$ is the identity map, then 
Theorem~\ref{tsesqui499} is a generalization of Corollary~3.9 in \cite{ERS4},
which followed from \cite{Kat1} Theorem~VIII.3.11.
Note that \cite{Kat1} Theorem~VIII.3.11 is a special case of 
Theorem~\ref{tsesqui499}.
The point in the following corollary is that the form $a$ is merely 
$j$-sectorial, but not necessarily $j$-elliptic.
It allows one to describe the associated operator also by a limit of
suitable perturbations.
This also underlines that the associated operator as we define it is
the natural object.

\begin{cor} \label{csesqui498}
Let $V,H$ be Hilbert spaces and $j \in \cl(V,H)$ with $j(V)$ dense in $H$.
Let $a \colon V \times V \to \Ci$ be a continuous $j$-sectorial form
with vertex~$\gamma$.
Let $b \colon V \times V \to \Ci$ be a $j$-elliptic continuous form.
Suppose that there exists a $\theta \in [0,\frac{\pi}{2})$
such that 
$b(u) \in \Sigma_\theta$
for all $u \in V$.
For all $n \in \Ni$ define $a_n = a + \frac{1}{n} \, b$.
Then $a_n$ is $j$-elliptic.
Let $A_n$ be the operator associated with $(a_n,j)$ and  
$A$ the operator associated with $(a,j)$.
Then 
\[
\lim_{n \to \infty} (\lambda I + A_n)^{-1} f= (\lambda I + A)^{-1} f
\]
for all $\lambda > (-\gamma) \vee 0$ and $f \in H$.
\end{cor}

We next consider the classical form associated with the $m$-sectorial operator~$A$.

\begin{prop} \label{psesqui361}
Let $a$ be a sesquilinear form, 
$H$ a Hilbert space and $j \colon D(a) \to H$ a linear map.
Suppose the form $a$ is $j$-sectorial and $j(D(a))$ is dense in $H$.
Let $A$ be the operator associated with $(a,j)$.
Then one has the following.
\begin{tabel}
\item \label{psesqui361-1}
There exists a unique closable sectorial form $a_r$ with form domain 
$j(D(a))$ such that $A$ is associated with $\overline{a_r}$.
\item \label{psesqui361-2}
$D(\overline{a_r}) = \{ x \in H : \mbox{there exists a bounded sequence }
   u_1,u_2,\ldots \mbox{ in } D(a) \mbox{ such that } x = \lim_{n \to \infty} j(u_n) \mbox{ in } H \} $.
\item \label{psesqui361-3}
There exists a $c > 0$ such that $\|j(u)\|_{a_r} \leq c \, \|u\|_a$ for all $u \in D(a)$.
In particular, if $D$ is a core for $a$, then $j(D)$ is a core for $\overline{a_r}$.
\item \label{psesqui361-4}
Let $h$ be the real part of $a$ and let $h_r$ be defined similarly as 
in Statement~{\rm \ref{psesqui361-1}}.
Then $D(\overline{a_r}) = D(\overline{h_r})$.
\end{tabel}
\end{prop}
\proof\
`\ref{psesqui361-1}'.
We use the notation as in the proof of Theorem~\ref{tsesqui340}.
Let $b$ be the closed sectorial form associated with $A$, i.e.\ the 
classical  form associated with $(\tilde a, \tilde j)$ given in 
Theorem~\ref{tsesqui212}\ref{tsesqui212-2}.
So $D(b) = \tilde j(V) = \tilde j(V(\tilde a))$ and 
$b(\tilde j(u),\tilde j(v)) = \tilde a(u,v)$ 
for all $u,v \in V(\tilde{a})$.
Then $j(D(a)) = \tilde j (q(D(a))) \subset \tilde j(V) = D(b)$.
We show that $j(D(a))$ is a core for $b$.
Let $x\in D(b)$. 
There exists a unique $u \in V(\tilde a)$
such that $\tilde j(u)=x$. 
There exist $u_1,u_2,\ldots \in D(a)$ such that 
$\lim q(u_n) = u$ in $V$. 
Let  $\pi_2$ be the projection of $V$ onto 
$V(\tilde{a})$ along the decomposition $V = \ker \tilde j \oplus V(\tilde{a})$.
Clearly $\pi_2(u) = u$.
In addition, $\pi_2$ is continuous and
$j(u_n) = \tilde j(q(u_n)) = \tilde j(\pi_2(q(u_n)))$ for all $n \in \Ni$.
Therefore
$\|x - j(u_n)\|_{D(b)}
= \|\pi_2(u) - \pi_2(q(u_n))\|_{V(\tilde{a})}
\leq \|\pi_2\| \, \|u - q(u_n)\|_V$
for all $n \in \Ni$, from which one deduces
that $\lim j(u_n) = x$ in $D(b)$.
We have shown that $j(D(a))$ is a core for $b$.
Let $a_r = b|_{j(D(a)) \times j(D(a))}$.
Then $b = \overline{a_r}$.
This proves existence of $a_r$. 
The uniqueness is obvious from \cite{Kat1} Theorem~VI.2.7.

`\ref{psesqui361-2}'.
`$\subset$'.
Let $x \in D(\overline{a_r}) = D(b)$.
Let $u_1,u_2,\ldots \in D(a)$ and $u \in V(\tilde a)$ be as in the 
proof of Statement~\ref{psesqui361-1}.
Then $\lim j(u_n) = x$ in $D(b)$, therefore also in $H$.
Moreover, $\lim q(u_n) = u$ in $V$.
So the sequence $q(u_1),q(u_2),\ldots$ is bounded in $V$.
But $\|u_n\|_a = \|q(u_n)\|_V$ for all $n \in \Ni$.
Thus the sequence $u_1,u_2,\ldots$ satisfies the requirements.

`$\supset$'.
Let $u_1,u_2,\ldots$ be a bounded sequence in $D(a)$, $x \in H$ and suppose that 
$\lim j(u_n) = x$ in $H$.
Then $q(u_1), q(u_2),\ldots$ is a bounded sequence in $V$.
So passing to a subsequence if necessary, there exists a 
$v \in V$ such that $\lim q(u_n) = v$ weakly in~$V$.
Then $\tilde j(v) = \lim j(u_n)$ weakly in~$H$.
Hence $x = \tilde j(v) \in \tilde j(V) = D(\overline{a_r})$.

`\ref{psesqui361-4}'.
The construction in the proof of Theorem~\ref{tsesqui340} with $h$ instead of $a$
leads to the same closed space $W$, then the same normed space $V_0$ and 
the same Banach space $V$.
Let $\tilde h \colon V \times V \to \Ci$ be the unique 
continuous form on $V$ such that $\tilde h(q(u), q(v)) = h(u,v)$ for all 
$u,v \in V$.
Then $\tilde h = \Re \tilde a$, the real part of $\tilde a$.
Let $h_c$ be the classical form associated with $(\tilde h,\tilde j)$.
Then $h_c = \overline{h_r}$ and $b = \overline{a_r}$ by 
part~\ref{psesqui361-1}.
Then Statement~\ref{psesqui361-4} follows from Proposition~\ref{psesqui214}.

`\ref{psesqui361-3}'.
Again by Proposition~\ref{psesqui214} there exists a $c \geq 1$ 
such that $\RRe b(x) \leq c \, h_c(x)$ for all $x \in \tilde j(V)$.
But $h_c(\tilde j(u)) \leq \tilde h(u) = \RRe \tilde a(u)$
for all $u \in V$.
So $\|\tilde j(u)\|_b \leq c \, \|u\|_{\tilde a}$ for all $u \in V$.
Then $\|j(u)\|_{a_r} \leq c \, \|q(u)\|_{\tilde a} = c \, \|u\|_a$
for all $u \in D(a)$.
The last assertion in Statement~\ref{psesqui361-3} is an 
immediate consequence.\hfill$\Box$

\ruimte

We call $a_r$ the {\bf regular part} and $\overline{a_r}$ 
the {\bf relaxed form} of the $j$-sectorial
form $a$.
If $D(a) \subset H$ and $j$ is the identity map, then it follows
from the proof of Proposition~\ref{psesqui361}\ref{psesqui361-1} that 
$a_r(x) = \tilde a(\pi_2(x))$ for all $x \in D(a)$, 
with the notation introduced there.
So if in addition $a$ is symmetric and {\bf positive}, i.e.\ if the numerical range 
$ \{a(u) : u \in D(a) \} $ is contained in $[0,\infty)$, then 
this terminology coincides with the one employed by 
Simon \cite{bSim5}, Section~2.
Under these assumptions
Simon characterized the regular part of $a$ as the 
largest closable form lying below $a$ for the order relation
$b_1 \leq b_2$ if and only if $D(b_2)\subset D(b_1)$
and $b_1(u)\leq b_2(u)$ for all $u\in D(b_2)$.
Of course, such an order relation does not exist
for sectorial forms. 
It seems to us, though, that the direct formula in 
Theorem \ref{tsesqui322}    expressing the generator 
directly in terms of the form $a$, is frequently more useful 
than the computation of $a_r$.
For positive $a$ Simon proved Proposition~\ref{psesqui361}\ref{psesqui361-2}
in \cite{bSim4}, Theorem~3.
Note that for general $a$ (but still $j$ the inclusion), 
the form $a$ is closable if and only if $a_r$ coincides
with $a$ on $D(a)$.

Let $a$ be a densely defined sectorial 
form and $A$ its associated operator, as above. 
If the form $a$ is symmetric, then
the associated operator $A$ is self-adjoint. 
But the converse is not true if the form $a$ is not closable. 
In order to see this, it suffices to consider the form
$(1+i)a$ where $a$ is the form as in Example~\ref{xsesqui323} below.

For general $j$-sectorial forms we also consider invariance of 
closed convex subsets.

\begin{prop} \label{psesqui330}
Let $a$ be a sesquilinear form, 
$H$ a Hilbert space and $j \colon D(a) \to H$ a linear map.
Suppose the form $a$ is accretive, $j$-sectorial and $j(D(a))$ is dense in $H$.
Let $A$ be the operator associated with $(a,j)$ and 
$S$ the semigroup generated by $-A$.
Let $C \subset H$ be a non-empty closed convex set and $P \colon H \to C$ the 
orthogonal projection.
Then the following are equivalent.
\begin{tabeleq}
\item \label{psesqui330-0.5}
$S_t C \subset C$ for all $t > 0$.
\item \label{psesqui330-1}
For all $u \in D(a)$ there exists a Cauchy sequence $w_1,w_2,\ldots$ in 
$(D(a),\|\cdot\|_a)$ such that 
$\lim_{n \to \infty} j(w_n) = P j(u)$ in $H$ and 
$\lim_{n \to \infty} \RRe a(w_n, u - w_n) \geq 0$.
\item \label{psesqui330-2}
For all $u \in D(a)$ there exists a bounded sequence $w_1,w_2,\ldots$ in 
$(D(a),\|\cdot\|_a)$ such that 
$\lim_{n \to \infty} j(w_n) = P j(u)$ in $H$ and 
$\limsup_{n \to \infty} \RRe a(w_n, u - w_n) \geq 0$.
\end{tabeleq}
\end{prop}
\proof\
We use the notation as in the proof of Theorem~\ref{tsesqui340}.
Clearly the form $\tilde a$ is accretive by continuity and density of~$V_0$.
We shall prove the equivalence with 
Condition~\ref{psesqui202-4} in Proposition~\ref{psesqui202}
for $D = V_0 = q(D(a))$, $\tilde a$ and $\tilde j$.

`\ref{psesqui330-0.5}$\Rightarrow$\ref{psesqui330-1}'.
Let $u \in D(a)$.
By Proposition~\ref{psesqui202}\ref{psesqui202-1}$\Rightarrow$\ref{psesqui202-4} 
there exists a $w \in V$ such that $\tilde j(w) = P j(u)$ 
and $\RRe \tilde a(w,q(u) - w) \geq 0$.
There are $w_1,w_2,\ldots \in D(a)$ such that $\lim q(w_n) = w$ in $V$.
Then the sequence $w_1,w_2,\ldots$ satisfies the requirements.

`\ref{psesqui330-1}$\Rightarrow$\ref{psesqui330-2}'.
Trivial.

`\ref{psesqui330-2}$\Rightarrow$\ref{psesqui330-0.5}'.
Let $u \in D(a)$.
By assumption there exists a bounded sequence
$w_1,w_2,\ldots$ in $(D(a), \|\cdot\|_a)$ such that 
$\lim j(w_n) = P j(u)$ in $H$ and 
$\limsup_{n \to \infty} a(w_n, u - w_n) \geq 0$.
Then $q(w_1),q(w_2),\ldots$ is a bounded sequence in $V$, so passing to a 
subsequence if necessary, it follows that it is weakly convergent.
Let $w = \lim_{n \to \infty} q(w_n)$ weakly in $V$.
Then $\tilde j(w) = \lim j(w_n)$ weakly in $H$, 
so $\tilde j(w) = P j(u) = P \tilde j(q(u))$.
Moreover, 
$\tilde a(w,q(u)) = \lim \tilde a(q(w_n), q(u))$
and $\RRe \tilde a(w,w) = \Re \tilde a(w) \leq \liminf \Re \tilde a(q(w_n))$
by \cite{Kat1} Lemma VIII.3.14a.
So $\RRe \tilde a(w,q(u) - w) \geq \limsup_{n \to \infty} \RRe a(w_n,u - w_n) \geq 0$.
Then Condition~\ref{psesqui330-0.5} follows from 
Proposition~\ref{psesqui202}\ref{psesqui202-4}$\Leftrightarrow$\ref{psesqui202-1}.\hfill$\Box$

\ruimte

\begin{remarkn} \label{rsesqui330.5}
Clearly Condition~\ref{psesqui330-1} in Proposition~\ref{psesqui330}
is valid if for all $u \in D(a)$ there exists a $w \in D(a)$
such that $j(w) = P j(u)$ and $\RRe a(w,u-w) \geq 0$.
\end{remarkn}

Proposition~\ref{psesqui330} has several consequences which will be useful
for differential operators in the next section.
If $(X,\cb,m)$ is a measure space and $a$ is a sesquilinear form in 
$L_2(X)$, then
we call $a$ {\bf real} if $\RRe u \in D(a)$ and $a(\RRe u) \in \Ri$
for all $u \in D(a)$.

\begin{cor} \label{csesqui331}
Let $(X,\cb,m)$ be a measure space and $a$ a 
densely defined sectorial form in $L_2(X)$.
Let $A$ be the operator associated with 
$a$ as in Theorem~{\rm \ref{tsesqui322}} and let 
$S$ be the semigroup generated by the operator $-A$.
\begin{tabel}
\item \label{csesqui331-1}
If $a$ is real, then $S_t L_2(X,\Ri) \subset L_2(X,\Ri)$ for all 
$t > 0$.
\item \label{csesqui331-2}
If $a$ is real, $u^+ \in D(a)$ and $a(u^+,u^-) \leq 0$ for all 
$u \in D(a) \cap L_2(X,\Ri)$, then $S$ is positive.
In particular, $|S_t u| \leq S_t |u|$ for all $t > 0$ and $u \in L_2(X)$.
\item \label{csesqui331-3}
If $a$ is accretive, real, $u \wedge \one \in D(a)$
and $a(u \wedge \one, (u-\one)^+) \geq 0$ for all $u \in D(a) \cap L_2(X,\Ri)$,
then $S$ is {\bf submarkovian}, i.e.,
$\|S_t u\|_\infty \leq \|u\|_\infty$ for all $u \in L_2(X) \cap L_\infty(X)$ and $t > 0$.
\item \label{csesqui331-4}
If $a$ is accretive, real, $u \wedge \one \in D(a)$ and
$a((u-\one)^+, u \wedge \one) \geq 0$ for all $u \in D(a) \cap L_2(X,\Ri)$,
then $\|S_t u\|_1 \leq \|u\|_1$ for all $u \in L_1(X) \cap L_2(X)$.
\end{tabel}
\end{cor}
\proof\
`\ref{csesqui331-1}'.
Replacing $a$ by $(u,v) \mapsto a(u,v) + \gamma \, (u,v)_H$
we may assume that $a$ is accretive.
Let $u \in D(a)$.
Set $w = \RRe u$.
Then $w \in D(a)$ and 
$\RRe a(w,u-w) = \RRe a(\RRe u, i \IIm u) = \IIm a(\RRe u, \IIm u) = 0$.
So by Proposition~\ref{psesqui330} the set $L_2(X,\Ri)$ is 
invariant under $S$.
(See also Remark~\ref{rsesqui330.5}.)

`\ref{csesqui331-2}'.
Again we may assume that $a$ is accretive.
Let $C = \{ u \in L_2(X,\Ri) : u \geq 0 \} $. 
Then $C$ is closed and convex.
Let $P$ be the orthogonal projection of $L_2(X)$ onto $C$.
Let $u \in D(a)$.
Then $Pu = (\RRe u)^+ \in D(a)$.
Moreover, 
$\RRe a(Pu, u - Pu)
= \RRe a((\RRe u)^+ , -(\RRe u)^- + i \IIm u)
= - a((\RRe u)^+, (\RRe u)^-) 
\geq 0$.
So by Proposition~\ref{psesqui330} the set $C$ is 
invariant under $S$.

`\ref{csesqui331-3}'. 
Let $C = \{ u \in L_2(X,\Ri) : u \leq \one \} $.
Then $C$ is closed and convex in $L_2(X)$.
The orthogonal projection $P \colon L_2(X) \to C$ is given by
$P u = (\RRe u) \wedge \one$.
It follows by assumption and Proposition~\ref{psesqui330} 
that the set $C$ is 
invariant under $S$.
Hence $S$ is submarkovian.

`\ref{csesqui331-4}'.
This follows by duality from Statement~\ref{csesqui331-3} and 
Remark~\ref{rsesqui346}.\hfill$\Box$

\ruimte

We end this section with an example which shows that in general (\ref{eSsesui2;1})
is restricted to $u \in D_H(\tilde a)$.

\begin{exam} \label{xsesqui323}
Let $H = L_2(0,1)$, $D(a) = C[0,1]$ and
\[
a(u,v) = \sum_{n=1}^\infty 2^{-n} \, u(q_n) \, \overline{v(q_n)}
\]
where $\{ q_n : n \in \Ni \} = [0,1] \cap \Qi$ with $q_n \neq q_m$ 
for all $n,m \in \Ni$ with $n \neq m$.
Moreover, let $j$ be the inclusion map.
Then it is not hard to characterize the completion of $D(a)$
and to show that the operator $A$ associated with~$a$ is the 
zero operator.\hfill$\Box$
\end{exam}

\section{Applications} \label{Ssesqui4}

We illustrate the theorems of the previous sections by several examples.

\subsection{Sectorial differential operators} \label{Ssesqui4.1}

First we consider differential operators on open sets in $\Ri^d$. 
We emphasize that the operators do not have to be symmetric
and may have complex coefficients.
The next lemma, whose proof is trivial, 
provides an efficient way to construct sectorial operators.

\begin{lemma} \label{lsesqui410}
Let $\Omega \subset \Ri^d$ be open.
For all $i,j \in \{ 1,\ldots,d \} $ let $a_{ij} \in L_{1,\loc}(\Omega)$.
Let $D(a)$ be a subspace of $L_2(\Omega)$ 
with $C_c^\infty(\Omega) \subset D(a)$.
Assume that $\partial_i u \in L_{1,\loc}(\Omega)$ as a distribution and
\[
\int_\Omega |(\partial_i u) \, a_{ij} \, \partial_j v|
< \infty
\]
for all $u,v \in D(a)$ and $i,j \in \{ 1,\ldots,d \} $.
Define the form $a \colon D(a) \times D(a) \to \Ci$ by
\[
a(u,v) = \sum_{i,j=1}^d \int_\Omega (\partial_i u) \, a_{ij} \, \overline{\partial_j v}
.  \]
Let $\theta \in [0,\frac{\pi}{2})$ and assume that 
$\sum_{i,j=1}^d a_{ij}(x) \, \xi_i \, \overline{\xi_j} \in \Sigma_\theta$
for all $\xi \in \Ci^d$ and a.e.\ $x \in \Omega$.
Then the form $a$ is sectorial with vertex $0$ and semi-angle $\theta$.
\end{lemma}

We call an operator $A$ associated with a form $a$
which satisfies the assumptions of Lemma~\ref{lsesqui410} a
{\bf sectorial differential operator} and $a$ a {\bf sectorial differential form}.
Then $-A$ generates a holomorphic semigroup.

The assumptions on the domain $D(a)$ and the coefficients
$a_{ij}$ are very general.
For example one can choose $D(a) = C_c^\infty(\Omega)$ together
with the condition $a_{ij} \in L_{1,\loc}(\Omega)$, or alternatively
if $a_{ij} \in L_\infty(\Omega)$ one can choose for 
$D(a)$ any subspace of 
$H^1(\Omega)$ with $C_c^\infty(\Omega) \subset D(a)$.

In order to avoid too many cases we will not consider unbounded
coefficients in this paper.
Let $\Omega \subset \Ri^d$ be open.
For all $i,j \in \{ 1,\ldots,d \} $ let $a_{ij} \in L_\infty(\Omega)$.
Define the form $a \colon D(a) \times D(a) \to \Ci$ by
\[
a(u,v) = \sum_{i,j=1}^d \int_\Omega (\partial_i u) \, a_{ij} \, \overline{\partial_j v}
,  \]
where $D(a)$ is a subspace of $H^1(\Omega)$ 
with $C_c^\infty(\Omega) \subset D(a)$.

We call $(a_{ij})$ {\bf strongly elliptic} if there exists a $\mu > 0$ such that 
$\RRe \sum_{i,j=1}^d a_{ij}(x) \, \xi_i \, \overline{\xi_j} \geq \mu \, |\xi|^2$
for all $\xi \in \Ci^d$ and a.e.\ $x \in \Omega$.
Clearly if $(a_{ij})$ is strongly elliptic, then there exists a
$\theta \in [0,\frac{\pi}{2})$ such that 
$\sum_{i,j=1}^d a_{ij}(x) \, \xi_i \, \overline{\xi_j} \in \Sigma_\theta$
for all $\xi \in \Ci^d$ and a.e.\ $x \in \Omega$.
We then also say that the form $a$ and associated operator are strongly elliptic.

Let $\theta \in [0,\frac{\pi}{2})$ and 
suppose that
$\sum_{i,j=1}^d a_{ij}(x) \, \xi_i \, \overline{\xi_j} \in \Sigma_\theta$
for all $\xi \in \Ci^d$ and a.e.\ $x \in \Omega$.
Let $l \colon D(a) \times D(a) \to \Ci$ be defined by 
$l(u,v) = \sum_{i=1}^d \int_\Omega \partial_i u \, \overline{\partial_i v}$.
For all $n \in \Ni$ let $a^{(n)} = a + \frac{1}{n} \, l$.
Although $a$ is not strongly elliptic in general,
the form $a^{(n)}$ is strongly elliptic for all $n \in \Ni$.
If $A$, $A_n$, $S$ and $S^{(n)}$ are the associated operators and semigroups,
then the conditions of Theorem~\ref{tsesqui499} are satisfied.
In particular the $A_n$ converge to $A$ strongly in the resolvent sense
and therefore $S^{(n)}_t$ converges strongly to $S_t$ for all $t > 0$.

We next show that under a mild condition on the form domain $D(a)$
the semigroup associated with a sectorial differential operator
satisfies Davies--Gaffney bounds.
If $F$ and $G$ are two non-empty subsets of $\Ri^d$, then 
$d(F,G)$ denotes the Euclidean distance.
The value of $M$ can be improved significantly if the coefficients are real.
(See \cite{ERSZ1} Proposition 3.1.)
In this paper the following version for complex coefficients suffices.

\begin{thm} \label{tsesqui411}
Let $\Omega \subset \Ri^d$ be open.
For all $i,j \in \{ 1,\ldots,d \} $ let $a_{ij} \in L_\infty(\Omega)$.
Let $\theta \in [0,\frac{\pi}{2})$.
Suppose 
$\sum_{i,j=1}^d a_{ij}(x) \, \xi_i \, \overline{\xi_j} \in \Sigma_\theta$
for all $\xi \in \Ci^d$ and a.e.\ $x \in \Omega$.
Define the form $a \colon D(a) \times D(a) \to \Ci$ by
\[
a(u,v) = \sum_{i,j=1}^d \int_\Omega (\partial_i u) \, a_{ij} \, \overline{\partial_j v}
,  \]
where $D(a)$ is a subspace of $H^1(\Omega)$ 
with $C_c^\infty(\Omega) \subset D(a)$.
Suppose $e^{\rho \, \psi}u \in D(a)$ for all $u \in D(a)$, $\rho \in \Ri$ and 
$\psi \in C_{\rm b}^\infty(\Ri^d,\Ri)$.
Let $S$ be the semigroup associated with $a$.
Then 
\begin{equation}
|(S_t u, v)|
\leq e^{-\frac{d(\Omega_1,\Omega_2)^2}{4 M t}} \, \|u\|_2 \, \|v\|_2
\label{etsesqui411;3}
\end{equation}
for all non-empty open $\Omega_1,\Omega_2 \subset \Omega$, 
$u \in L_2(\Omega_1)$, $v \in L_2(\Omega_2)$ and $t > 0$, 
where $M = 3(1 + \tan \theta)^2 (1 + \sum_{i,j=1}^d \|a_{ij}\|_\infty)$.
\end{thm}
\proof\
First suppose that $(a_{ij})$ is strongly elliptic.
Let $\rho > 0$ and $\psi \in C_{\rm b}^\infty(\Ri^d,\Ri)$
with $\|\nabla \psi\|_\infty \leq 1$.
Define the form $a_\rho \colon D(a) \times D(a) \to \Ci$ by
\[
a_\rho(u,v) = \sum_{i,j=1}^d \int_\Omega (\partial_i u + \rho \, \psi_i \, u) 
             \, a_{ij} \, \overline{\partial_j v - \rho \, \psi_j \, v}
,  \]
where $\psi_i = \partial_i \psi$ for all $i \in \{ 1,\ldots,d \} $.
Then 
\begin{eqnarray}
\RRe a_\rho(u)
& = & \RRe a(u)
   + \rho \, \RRe \int_\Omega \sum_{i,j=1}^d \psi_i \, u \, a_{ij} \, \overline{\partial_j u}
   - \rho \, \RRe \int_\Omega \sum_{i,j=1}^d (\partial_i u) \, a_{ij} \, \psi_j \, \overline{u}  \nonumber  \\
& & \hspace{10mm} {}
   - \rho^2 \, \RRe \int_\Omega \sum_{i,j=1}^d \psi_i \, a_{ij} \, \psi_j \, |u|^2
\label{etsesqui411;1}
\end{eqnarray}
for all $u \in D(a)$.
It follows from the estimate (1.15) of Subsection~VI.1.2 in \cite{Kat1} that 
\begin{eqnarray*}
\Big| \sum_{i,j=1}^d a_{ij}(x) \, \xi_i \, \overline{\eta_j} \Big|
& \leq & (1 + \tan \theta) \Big( \RRe \sum_{i,j=1}^d a_{ij}(x) \, \xi_i \, \overline{\xi_j} \Big)^{1/2}
     \Big( \RRe \sum_{i,j=1}^d a_{ij}(x) \, \eta_i \, \overline{\eta_j} \Big)^{1/2}  \\
& \leq & \varepsilon \RRe \sum_{i,j=1}^d a_{ij}(x) \, \xi_i \, \overline{\xi_j} 
   + \frac{(1 + \tan \theta)^2}{4 \varepsilon} \, 
          \RRe \sum_{i,j=1}^d a_{ij}(x) \, \eta_i \, \overline{\eta_j}
\end{eqnarray*}
for all $\xi,\eta \in \Ci^d$, $\varepsilon > 0$ and a.e.\ $x \in \Omega$.
Choosing $\xi_i = (\partial_i u)(x)$, $\eta_i = (\psi_i \, u)(x)$ and 
$\varepsilon = \frac{1}{4 \rho}$ it follows that 
\[
\rho \Big| \int_\Omega \sum_{i,j=1}^d (\partial_i u) \, a_{ij} \, \psi_j \, \overline{u} \Big|
\leq \tfrac{1}{4} \RRe a(u) 
   +  (1 + \tan \theta)^2 \, \rho^2 \, 
     \RRe \int_\Omega \sum_{i,j=1}^d \psi_i \, a_{ij} \, \psi_j \, |u|^2
.  \]
Similarly the second term in (\ref{etsesqui411;1}) can be estimated.
Hence 
\begin{eqnarray}
\RRe a_\rho(u) 
& \geq & \tfrac{1}{2} \RRe a(u) 
   - \Big( 1 + 2 (1 + \tan \theta)^2 \Big) \rho^2 \, 
     \RRe \int_\Omega \sum_{i,j=1}^d \psi_i \, a_{ij} \, \psi_j \, |u|^2  \nonumber  \\
& \geq & \tfrac{1}{2} \RRe a(u) - M \, \rho^2 \, \|u\|_2^2
.
\label{etsesqui411;2}
\end{eqnarray}
Define $U_{\pm \rho} \colon L_2(\Omega) \to L_2(\Omega)$ by
$U_{\pm \rho} v = e^{\pm \rho \, \psi} v$.
Then $U_{\pm \rho} D(a) \subset D(a)$.
Moreover, $a_\rho(u,v) = a(U_\rho u, U_{-\rho} v)$ for all $u,v \in D(a)$.
Since $(a_{ij})$ is strongly elliptic, the forms $a$ and $a_\rho$ are 
sectorial (cf.\ Lemmas~3.6 and 3.7 in \cite{AE1}).
Let $A$ and $A_\rho$ be the associated operators and let 
$S^{(\rho)}$ be the semigroup generated by $-A_\rho$.
Then $A_\rho = U_{-\rho} \, A \, U_\rho$ and $S^{(\rho)}_t = U_{-\rho} \, S_t \, U_\rho$
for all $t > 0$.
It follows from (\ref{etsesqui411;2}) that 
\begin{equation}
\|S^{(\rho)}_t\|_{2 \to 2} \leq e^{M \, \rho^2 \, t}
\label{etsesqui411;5}
\end{equation}
for all $t > 0$.
Then 
\begin{eqnarray*}
|(S_t u, v)|
& = & |(S^{(\rho)}_t \, U_{-\rho} u , U_\rho v)|
\leq \|S^{(\rho)}_t\|_{2 \to 2} \, \|U_{-\rho} u\|_2 \, \|U_\rho v\|_2
\leq e^{M \, \rho^2 \, t} \, e^{-\rho \, d_\psi(\Omega_1,\Omega_2)} \, \|u\|_2 \, \|v\|_2
\end{eqnarray*}
for all $u \in L_2(\Omega_1)$ and $v \in L_2(\Omega_2)$,
where $d_\psi(\Omega_1,\Omega_2) = \inf_{x \in \Omega_1} \psi(x) - \sup_{x \in \Omega_2} \psi(x)$.
Minimizing over all $\psi \in C_{\rm b}^\infty(\Ri^d)$ with $\|\nabla \psi\|_\infty \leq 1$
one deduces that 
\[
|(S_t u, v)|
\leq e^{M \, \rho^2 \, t} \, e^{-\rho \, d(\Omega_1,\Omega_2)} \, \|u\|_2 \, \|v\|_2
\]
and choosing $\rho = \frac{d(\Omega_1,\Omega_2)}{2Mt}$ gives
\[
|(S_t u, v)|
\leq e^{- \frac{d(\Omega_1,\Omega_2)^2}{4 M t}} \, \|u\|_2 \, \|v\|_2
\]
for all $u \in L_2(\Omega_1)$, $v \in L_2(\Omega_2)$ and $t > 0$.

Finally we drop the assumption that $(a_{ij})$ is strongly elliptic.
For all $n \in \Ni$ define $a^{(n)}_{ij} = a_{ij} + \frac{1}{n} \, \delta_{ij}$.
Then $(a^{(n)}_{ij})$ is strongly elliptic.
If $S^{(n)}$ is the associated semigroup, then 
$\lim_{n \to \infty} S^{(n)}_t = S_t$ strongly for all $t > 0$ by 
Theorem~\ref{tsesqui499}.
Hence the theorem follows.\hfill$\Box$

\ruimte

We next consider locality properties of the relaxed form 
$\overline{a_r}$ of the 
sectorial form $a$.

\begin{cor} \label{csesqui413}
Assume the notation and assumptions of Theorem~{\rm \ref{tsesqui411}}.
Then $\overline{a_r}(u,v) = 0$ for all $u,v \in D(\overline{a_r})$ 
with disjoint compact supports.
\end{cor}
\proof\
There exist open non-empty $\Omega_1,\Omega_2 \subset \Ri^d$ 
such that $\supp u \subset \Omega_1$, $\supp v \subset \Omega_2$ 
and $d(\Omega_1,\Omega_2) > 0$.
Then it follows from Theorem~\ref{tsesqui411} that there exists a 
$b > 0$ such that 
\[
|((I - S_t)u,v)|
= |(S_t u, v)|
\leq e^{-b t^{-1}} \|u\|_2 \, \|v\|_2
\]
for all $t > 0$.
Hence by \cite{Ouh5} Lemma~1.56 one deduces that 
\[
|\overline{a_r}(u,v)|
= \lim_{t \downarrow 0} t^{-1} |((I - S_t)u,v)|
\leq \lim_{t \downarrow 0} t^{-1} e^{-b t^{-1}} \|u\|_2 \, \|v\|_2
= 0
\]
as required.\hfill$\Box$

\ruimte

If $\Omega \subset \Ri^d$ define 
\[
L_{2,c}(\Omega)
= \{ u \in L_2(\Omega) : \supp u \mbox{ is compact} \}
.  \]
Another corollary of Theorem~{\rm \ref{tsesqui411}} is that 
$S_t$ maps $L_{2,c}(\Omega)$ into $L_1(\Omega)$.
This is a special case of the following lemma.

For all $R > 0$ let $B_R$ denote the open ball
in $\Ri^d$ with centre $0$ and radius~$R$.
Set $\chi_R = \one_{B_R}$.

\begin{lemma} \label{lsesqui414}
Let $d \in \Ni$.
There exists a constant $c_d > 0$ such that the following holds.
Let $\Omega \subset \Ri^d$ be open and 
$T \in \cl(L_2(\Omega))$.
Let $N > 0$ and suppose that 
\[
|(T u, v)|
\leq e^{-\frac{d(\Omega_1,\Omega_2)^2}{N}} \, \|u\|_2 \, \|v\|_2
\]
for all non-empty open $\Omega_1,\Omega_2 \subset \Omega$, 
$u \in L_2(\Omega_1)$ and $v \in L_2(\Omega_2)$.
Then $T L_{2,c}(\Omega) \subset L_1(\Omega)$ and
\[
\|(\one - \chi_{2R}) T u\|_1
\leq c_d \, R^{-1} N^{\frac{d+2}{4}} \,
    e^{- \frac{R^2}{2N}} \, \|u\|_2
\]
for all $R > 0$ and $u \in L_2(\Omega)$ with 
$\supp u \subset B_R$.
\end{lemma}
\proof\
Since $\chi_{2R} \, T u \in L_2(\Omega \cap B_{2R}) \subset L_1(\Omega)$
it suffices to show the estimate.
Let $\varphi \in C_c(\Omega)$.
Then
\begin{eqnarray*}
|((\one - \chi_{2R}) T u,\varphi)|
& = & |(T u, (\one - \chi_{2R}) \varphi)|  \\
& \leq & \sum_{n=1}^\infty |(T u, (\chi_{(n+2)R} - \chi_{(n+1)R}) \varphi)|  \\
& \leq & \sum_{n=1}^\infty e^{-\frac{n^2 R^2}{N}} 
    \, \|u\|_2 \, \|(\chi_{(n+2)R} - \chi_{(n+1)R}) \varphi\|_2  \\
& \leq & \sum_{n=1}^\infty e^{-\frac{n^2 R^2}{N}} 
    \, ((n+2) R)^{d/2} \, |B_1|^{1/2} \, \|u\|_2 \, \|\varphi\|_\infty  \\
& \leq & 3^{d/2} |B_1|^{1/2} \, e^{-\frac{R^2}{2N}} \|u\|_2 \, \|\varphi\|_\infty
    \sum_{n=1}^\infty e^{-\frac{n^2 R^2}{2N}} \, (n R)^{d/2}  
.
\end{eqnarray*}
Let $c' > 0$ be such that $x^{d/4} \leq c' e^x$ uniformly for all $x > 0$.
Then $c'$ can be chosen to depend only on $d$.
Note that 
$\sum_{n=1}^\infty e^{-a n^2} 
   \leq \int_0^\infty e^{-a x^2} \, dx 
   = \sqrt{\frac{\pi}{4a}}$
for all $a > 0$.
Therefore 
\begin{eqnarray*}
\sum_{n=1}^\infty e^{-\frac{n^2 R^2}{2N}} \, (n R)^{d/2} 
& = & (4N)^{d/4} \sum_{n=1}^\infty 
     e^{-\frac{n^2 R^2}{2N}} \, \Big( \frac{n^2 \, R^2}{4N} \Big)^{d/4}  \\
& \leq & c' \, (4N)^{d/4} \sum_{n=1}^\infty e^{-\frac{n^2 R^2}{4N}}
\leq c' \, (4N)^{d/4} \Big( \frac{\pi \, N}{R^2} \Big)^{1/2}
.  
\end{eqnarray*}
Then the lemma follows by taking the supremum 
over all $\varphi$ with $\|\varphi\|_\infty \leq 1$.\hfill$\Box$

\ruimte

As a consequence one deduces $L_1$-convergence of the approximate semigroups
on $L_{2,c}(\Omega)$.
Recall that the coefficients in Theorem~\ref{tsesqui411} are complex.

\begin{lemma} \label{lsesqui413.2}
Assume the notation and assumptions of Theorem~{\rm \ref{tsesqui411}}.
For all $n \in \Ni$ let $a^{(n)} = a + \frac{1}{n} \, l$, where 
$l$ is the form with $D(l) = D(a)$ and 
$l(u,v) = \sum_{i=1}^d \int_\Omega \partial_i u \, \overline{\partial_i v}$.
Let $S^{(n)}$ be the semigroup associated with $a^{(n)}$.
Then $\lim_{n \to \infty} S^{(n)}_t u = S_t u$ in $L_1(\Omega)$ for all $t > 0$
and $u \in L_{2,c}(\Omega)$.
\end{lemma}
\proof\
It follows from Theorem~{\rm \ref{tsesqui411}} that there exists an $M > 0$ such that 
\[
|(S_t u,v)| \vee |(S^{(n)}_t u,v)|
\leq e^{-\frac{d(\Omega_1,\Omega_2)^2}{4 M t}} \, \|u\|_2 \, \|v\|_2
\]
for all $n \in \Ni$, non-empty open $\Omega_1,\Omega_2 \subset \Omega$, 
$u \in L_2(\Omega_1)$, $v \in L_2(\Omega_2)$ and $t > 0$, 
Let $c_d > 0$ be as in Lemma~\ref{lsesqui414}.
Let $u \in L_{2,c}(\Omega)$ and $t > 0$.
Then 
\[
\|(\one - \chi_{2R}) S^{(n)}_t u\|_1
\leq c_d \, R^{-1} (4Mt)^{\frac{d+2}{4}} \,
    e^{- \frac{R^2}{8Mt}} \, \|u\|_2
\]
for all $n \in \Ni$ and $R > 0$ with 
$\supp u \subset B_R$.
So $\lim_{R \to \infty} (\one - \chi_{2R}) S^{(n)}_t u = 0$ in $L_1(\Omega)$
uniformly in $n \in \Ni$.
Similarly, $\lim_{R \to \infty} (\one - \chi_{2R}) S_t u = 0$ in $L_1(\Omega)$.
So it suffices to prove that $\lim_{n \to \infty} \chi_{2R} (S^{(n)}_t u - S_t u) = 0$
for large $R > 0$.
Since $\|\chi_{2R} (S^{(n)}_t u - S_t u)\|_1 \leq |B_{2R}|^{1/2} \, \|S^{(n)}_t u - S_t u\|_2$
for all $n \in \Ni$ and $R > 0$, it follows from Theorem~\ref{tsesqui499}
that $\lim_{n \to \infty} \chi_{2R} (S^{(n)}_t u - S_t u) = 0$ in 
$L_1(\Omega)$ for all $R > 0$.\hfill$\Box$

\ruimte

For strongly elliptic operators one can strengthen the conclusions
of Theorem~\ref{tsesqui411}.

\begin{lemma} \label{lsesqui413.3}
Assume the notation and assumptions of Theorem~{\rm \ref{tsesqui411}}.
In addition suppose that the operator is strongly elliptic.
Then one has the following.
\begin{tabel}
\item \label{lsesqui413.3-1}
$S_t L_2(\Omega) \subset H^1(\Omega)$ for all $t > 0$.
\item \label{lsesqui413.3-2}
There exist $c,M' > 0$ such that 
\[
|(\partial_i S_t u, v)|
\leq c \, e^{-\frac{d(\Omega_1,\Omega_2)^2}{M' t}} \, \|u\|_2 \, \|v\|_2
\]
for all non-empty open $\Omega_1,\Omega_2 \subset \Omega$, 
$u \in L_2(\Omega_1)$, $v \in L_2(\Omega_2)$ and $t > 0$.
\item \label{lsesqui413.3-3}
If $u \in L_{2,c}(\Omega)$, then 
$S_t u, \partial_i S_t u \in L_1(\Omega)$ for all 
$t > 0$ and $i \in \{ 1,\ldots,d \} $.
Moreover, $t \mapsto \|\partial_i S_t u\|_1$ is locally bounded.
\end{tabel}
\end{lemma}
\proof\
`\ref{lsesqui413.3-1}'.
Let $b$ be the sectorial differential form with form domain $D(b) = H^1(\Omega)$
and coefficients $a_{ij}$.
Since $(a_{ij})$ is strongly elliptic if follows that $b$ is closed.
Clearly $b$ is an extension of $a$.
So $a$ is closable.
Let $A$ be the operator associated with $a$.
Then $A$ is the operator associated with $\overline a$.
Since $S$ is holomorphic one deduces that 
$S_t L_2(\Omega) 
\subset D(A)
\subset D(\overline a)
\subset D(b)
= H^1(\Omega)$ for all $t > 0$.

`\ref{lsesqui413.3-3}'.
This is a consequence of Lemma~\ref{lsesqui414}
and the estimates of Theorem~\ref{tsesqui411} and Statement~\ref{lsesqui413.3-2}.

`\ref{lsesqui413.3-2}'.
We use the notation as in the proof of Theorem~\ref{tsesqui411}.
Fix $\theta' \in (\theta,\frac{\pi}{2})$.
For all $\varphi \in \Ri$ with $|\varphi| < \theta' - \theta$ 
define $a^{[\varphi]}_{ij} = e^{i \varphi} \, a_{ij}$ for all
$i,j \in \{ 1,\ldots,d \} $.
Then $\sum_{i,j=1}^d a^{[\varphi]}_{ij}(x) \, \xi_i \, \overline{\xi_j} \in \Sigma_{\theta'}$
for all $\xi \in \Ci^d$ and a.e.\ $x \in \Omega$.
Let $a^{[\varphi]}$ be the corresponding form with form domain $D(a)$.
For all $\rho > 0$ let $a^{[\varphi]}_\rho$, $A^{[\varphi]}$, $A^{[\varphi]}_\rho$,
$S^{[\varphi]}$ and $S^{[\varphi] \, \rho}$ be the form, operators and semigroups
defined naturally as in the proof of Theorem~\ref{tsesqui411}.
Then it follows from (\ref{etsesqui411;5}) that 
\[
\|S^{[\varphi] \, \rho}_t\|_{2 \to 2} 
\leq e^{M_1 \, \rho^2 \, t}
\]
for all $\rho,t > 0$ and $|\varphi| < \theta' - \theta$, where 
$M_1 = 3(1 + \tan \theta')^2 (1 + \sum_{i,j=1}^d \|a_{ij}\|_\infty)$.
But $S^{[\varphi] \, \rho}_t = e^{-t e^{i \varphi} A_\rho} = S^\rho_{t e^{i \varphi}}$.
So $\|S^\rho_{t e^{i \varphi}}\|_{2 \to 2} \leq e^{M_1 \, \rho^2 \, t}$
for all $t,\rho > 0$ and $|\varphi| < \theta' - \theta$.
Since $S^\rho$ is a holomorphic semigroup on the interior of 
$\Sigma_{\frac{\pi}{2} - \theta'}$ it follows that
\[
S^\rho_t
= \frac{1}{2 \pi i} \, \int_{\Gamma_r(t)} \frac{1}{z-t} \, S^\rho_z \, dz
\]
for all $t > 0$, where $\Gamma_r(t)$ is the circle centred at $t$ and radius 
$r = c \, t$ and $c = \sin \frac{1}{2} (\frac{\pi}{2} - \theta')$.
Therefore
\[
\|A_\rho S^\rho_t\|_{2 \to 2}
\leq \frac{1}{2 \pi}  \, \int_{\Gamma_r(t)} \frac{1}{|z-t|^2} \, \|S^\rho_z\|_{2 \to 2} \, d|z|
\leq \frac{1}{c \, t} \, e^{M_2 \rho^2 t}
\]
for all $\rho,t > 0$, where $M_2 = M_1(1+c)$.
It then follows from (\ref{etsesqui411;2}) that 
\begin{eqnarray*}
\tfrac{1}{2} \, \mu \sum_{i=1}^d \|\partial_i \, S^\rho_t u\|_2^2
& \leq & \RRe a_\rho(S^\rho_t u, S^\rho_t u) + M \rho^2 \, \|S^\rho_t u\|_2^2  \\
& \leq & \|A_\rho \, S^\rho_t u\|_2 \, \|S^\rho_t u\|_2 
      + M \rho^2 \, \|S^\rho_t u\|_2^2  \\
& \leq & \frac{1}{c \, t} \, e^{M_2 \, \rho^2 \, t} \, e^{M \, \rho^2 \, t} \|u\|_2^2
      + M \rho^2 \, e^{2 M \, \rho^2 \, t} \|u\|_2^2
.
\end{eqnarray*}
Hence there exist $c_3,M_3 > 0$ such that 
\[
\|\partial_i \, S^\rho_t\|_{2 \to 2} 
\leq c_3 \, t^{-1/2} \, e^{M_3 \, \rho^2 \, t}
\]
for all $i \in \{ 1,\ldots,d \} $ and $\rho,t > 0$.
Since 
\[
\|U_{-\rho} \, \partial_i \, S_t \, U_\rho\|_{2 \to 2}
= \|(\partial_i + \rho \, \psi_i) S^\rho_t\|_{2 \to 2}
\leq \|\partial_i \, S^\rho_t\|_{2 \to 2}
   + |\rho \, \psi_i| \, \|S^\rho_t\|_{2 \to 2}
\leq c_4 \, t^{-1/2} \, e^{M_4 \, \rho^2 \, t}
\]
for suitable $c_4,M_4 > 0$,
Statement~\ref{lsesqui413.3-2} follows as at the end of the proof
of Theorem~\ref{tsesqui411}.\hfill$\Box$

\ruimte

The conditions on the form domain in Theorem~{\rm \ref{tsesqui411}}
are satisfied in case of Neumann boundary conditions, i.e.\
if $D(a) = H^1(\Omega)$. 
We next show that if $D(a) = H^1(\Omega)$, then a strong locality property is valid.
We start with a lemma for (complex) strongly elliptic operators.

\begin{lemma} \label{lsesqui413.5}
Let $\Omega \subset \Ri^d$ be open.
For all $i,j \in \{ 1,\ldots,d \} $ let $a_{ij} \in L_\infty(\Omega)$.
Suppose $(a_{ij})$ is strongly elliptic.
Define $a \colon H^1(\Omega) \times H^1(\Omega) \to \Ci$ by
\[
a(u,v) = \sum_{i,j=1}^d \int_\Omega (\partial_i u) \, a_{ij} \, \overline{\partial_j v}
.  \]
Let $S$ be the semigroup associated with $a$.
Then $(S_t u,\one) = (u,\one)$ for all $u \in L_{2,c}(\Omega)$ and $t > 0$.
\end{lemma}
\proof\
Fix $\tau \in C_c^\infty(\Ri^d)$ such that $\tau|_{B_1} = \one$.
For all $n \in \Ni$ define $\tau_n \in C_c^\infty(\Ri^d)$ by
$\tau_n(x) = \tau(n^{-1} x)$.
For all $n \in \Ni$ define $f_n \colon (0,\infty) \to \Ci$ by
$f_n(t) = (S_t u, \tau_n \, \one_\Omega)$.
Note that $\tau_n \, \one_\Omega \in H^1(\Omega) = D(a)$ for all $n \in \Ni$.
Therefore 
\[
f_n'(t)
= - a(S_t u, \tau_n \, \one_\Omega)  
= - \sum_{i,j=1}^d (\partial_i \, S_t u, a_{ij} \, \partial_j(\tau_n \, \one_\Omega)) 
= - \sum_{i,j=1}^d (\partial_i \, S_t u, a_{ij} \, (\partial_j \tau_n) \, \one_\Omega)  
\]
and 
\[
|f_n'(t)|
\leq \sum_{i,j=1}^d \|\partial_i \, S_t u\|_1 \, \|a_{ij}\|_\infty \,
      n^{-1} \, \|\partial_j \tau\|_\infty
\]
for all $n \in \Ni$ and $t > 0$, where we used that 
$\partial_i \, S_t u \in L_1(\Omega)$ by Lemma~\ref{lsesqui413.3}\ref{lsesqui413.3-3}.
So $\lim_{n \to \infty} f_n'(t) = 0$ locally uniform on $(0,\infty)$.
In addition, $\lim_{n \to \infty} f_n(t) = (S_t u, \one)$ for all $t \in (0,\infty)$.
Therefore $t \mapsto (S_t u, \one)$ is constant.
Since $\lim_{t \downarrow 0} (S_t u,\one) = (u,\one)$ the lemma follows.\hfill$\Box$

\ruimte

We are now able to prove strong locality for Neumann sectorial differential operators.
Note that our conditions allow that the coefficients are $0$ 
on part or even the entire domain.

\begin{prop} \label{psesqui415}
Let $\Omega \subset \Ri^d$ be open.
For all $i,j \in \{ 1,\ldots,d \} $ let $a_{ij} \in L_\infty(\Omega)$.
Let $\theta \in [0,\frac{\pi}{2})$.
Suppose 
$\sum_{i,j=1}^d a_{ij}(x) \, \xi_i \, \overline{\xi_j} \in \Sigma_\theta$
for all $\xi \in \Ci^d$ and a.e.\ $x \in \Omega$.
Define the form $a$ with form domain $D(a) = H^1(\Omega)$ by
\[
a(u,v) = \sum_{i,j=1}^d \int_{\Ri^d} (\partial_i u) \, a_{ij} \, \overline{\partial_j v}
.  \]
Then one has the following.
\begin{tabel}
\item \label{psesqui415-1}
$\overline{a_r}(u,v) = 0$ for all $u,v \in D(\overline{a_r})$ 
with compact support such that $v$ is constant on a 
neighbourhood of the support of $u$.
\item \label{psesqui415-2}
If $S$ is the semigroup associated with $a$, then $(S_t u,\one) = (u,\one)$ 
for all $t > 0$ and $u \in L_{2,c}(\Omega)$.
\end{tabel}
\end{prop}
\proof\
We first prove Statement~\ref{psesqui415-2}.
For all $n \in \Ni$ let $a^{(n)} = a + \frac{1}{n} \, l$, where 
$l$ is the form with $D(l) = H^1(\Omega)$ and 
$l(u,v) = \sum_{i=1}^d \int_\Omega \partial_i u \, \overline{\partial_i v}$.
Let $S^{(n)}$ be the semigroup associated with $a^{(n)}$.
Then $(S_t u,\one) = \lim_{n \to \infty} (S^{(n)}_t u,\one) = (u,\one)$
for all $t > 0$ and $u \in L_{2,c}(\Omega)$
by Lemmas~\ref{lsesqui413.2} and \ref{lsesqui413.5}.

Next let $u,v \in D(\overline{a_r})$ 
with compact support such that $v$ is constant on a 
neighbourhood of the support of $u$.
Then there exist an open set $U$ and a $\lambda \in \Ci$ such 
that $\supp u \subset U$ and $v(x) = \lambda$ for all $x \in U$.
Therefore $(u,v) = \lambda \, (u,\one) = \lambda \, (S_t u, \one)$
for all $t > 0$.

Let $c_d > 0$ be the constant in Lemma~\ref{lsesqui414}, which 
depends only on $d$.
Moreover, set 
\[
M = 3(1 + \tan \theta)^2 (1 + \sum_{i,j=1}^d \|a_{ij}\|_\infty)
.  \]
Fix $R > 0$ such that $\supp u \subset B_R$.

Now let $t > 0$.
Then 
\begin{eqnarray*}
((I - S_t) u, v)
& = & \lambda \, (S_t u, \one) - (S_t u, v)  \\
& = & \lambda \, (S_t u, \one - \chi_{2R})
   + (S_t u, \lambda \, \chi_{2R} - v) 
.
\end{eqnarray*}
We estimate the terms separately.
First, $S$ satisfies  the 
Davies--Gaffney bounds (\ref{etsesqui411;3}) of Theorem~\ref{tsesqui411}.
So one estimates
\[
|(S_t u, \one - \chi_{2R})|
\leq \|(\one - \chi_{2R}) S_t u\|_1
\leq c_d \, R^{-1} (4 M \, t)^{\frac{d+2}{4}} \,
    e^{- \frac{R^2}{8 M t}} \, \|u\|_2
\]
by Lemma~\ref{lsesqui414}.
Next, let $D > 0$ be the distance between $\supp u$ and $U^{\rm c}$.
Then it follows from Theorem~\ref{tsesqui411} that 
\begin{eqnarray*}
|(S_t u, \lambda \, \chi_{2R} - v)|
& \leq & e^{-\frac{D^2}{4 M t}} \, \|u\|_2 \, \|\lambda \, \chi_{2R} - v\|_2  \\
& \leq & (|\lambda| \, (2R)^{d/2} \, |B_1|^{1/2} + \|v\|_2) 
     \, e^{-\frac{D^2}{4 M t}} \, \|u\|_2
.  
\end{eqnarray*}
Therefore 
\begin{eqnarray*}
t^{-1} |((I - S_t) u, v)|
& \leq & |\lambda| \, c_d \, R^{-1} t^{-1} (4 M \, t)^{\frac{d+2}{4}} \,
    e^{- \frac{R^2}{8 M t}} \, \|u\|_2  \\*
& & \hspace{20mm} {}
   + (|\lambda| \, (2R)^{d/2} \, |B_1|^{1/2} + \|v\|_2) 
     \, t^{-1} e^{-\frac{D^2}{4 M t}} \, \|u\|_2  
\end{eqnarray*}
for all $t > 0$.
Since 
$\overline{a_r}(u,v) = \lim_{t \downarrow 0} t^{-1} ((I - S_t)u,v)$
the proposition follows.\hfill$\Box$

\ruimte

Up to now the coefficients were allowed to be complex in this section.
If the coefficients are real, but possibly not symmetric,
then one has the following application of Corollary~\ref{csesqui331} and 
Proposition~\ref{psesqui415}.

\begin{cor} \label{csesqui412}
Let $\Omega \subset \Ri^d$ be open.
For all $i,j \in \{ 1,\ldots,d \} $ let $a_{ij} \in L_\infty(\Omega,\Ri)$.
Let $\theta \in [0,\frac{\pi}{2})$.
Suppose 
$\sum_{i,j=1}^d a_{ij}(x) \, \xi_i \, \overline{\xi_j} \in \Sigma_\theta$
for all $\xi \in \Ci^d$ and a.e.\ $x \in \Omega$.
Define the form $a \colon D(a) \times D(a) \to \Ci$ by
\[
a(u,v) = \sum_{i,j=1}^d \int_\Omega (\partial_i u) \, a_{ij} \, \overline{\partial_j v}
,  \]
where $D(a) = H^1(\Omega)$ or $D(a) = H^1_0(\Omega)$.
Let $S$ be the semigroup associated with $a$.
Then $S$ is real, positive and $S$ extends consistently to a 
contraction semigroup
on $L_p(\Omega)$ for all $p \in [1,\infty]$,
which is a $C_0$-semigroup if $p \in [1,\infty)$ and the 
adjoint of a $C_0$-semigroup if $p = \infty$.
Moreover, if $D(a) = H^1(\Omega)$, then $S_t \one_\Omega = \one_\Omega$ for all $t > 0$.
\end{cor}
\proof\
Only the last statement needs comments. 
Since $L_{2,c}(\Omega)$ is dense in $L_1(\Omega)$ one deduces from 
Proposition~\ref{psesqui415}\ref{psesqui415-2} 
that $(S_t u,\one) = (u,\one)$ for all $u \in L_1(\Omega)$.
Then the claim follows by duality and Remark~\ref{rsesqui346}.\hfill$\Box$

\ruimte

Thus for real coefficients and Neumann boundary conditions the semigroup
$S$ is stochastic on $L_1$.

\subsection{Multiplicative perturbation} \label{Ssesqui4.4}

We  perturb the Dirichlet Laplacian by choosing a special function~$j$.
Let $\Omega \subset \Ri^d$ be open and bounded.
Then we obtain a possibly degenerate operator as follows.

\begin{prop} \label{psesqui465}
Let $m \colon \Omega \to (0,\infty)$ be such that $\frac{1}{m} \in L_{2,\loc}(\Omega)$.
Define the operator, formally denoted by $(m \Delta m)$ on $L_2(\Omega)$
by the following.
Let $w,f \in L_2(\Omega)$.
Then we define $w \in D((m \Delta m))$ and $(m \Delta m) \, w = f$ if and only if
$m w \in H^1_0(\Omega)$ and $\Delta (m \, w) = \frac{f}{m}$ in $\cd(\Omega)'$.

Then the operator $(m \Delta m)$ is self-adjoint and $(m \Delta m)$ generates
a positive semigroup~$S$.
Moreover, the set 
\[
C = \{ f \in L_2(\Omega,\Ri) : f \leq \tfrac{1}{m} \} 
\]
is invariant under $S$.
\end{prop}
\proof\
Let $V = H^1_0(\Omega) \cap L_2(\Omega,\frac{1}{m^2} \, dx)$
and define $j \in \cl(V, L_2(\Omega))$ by 
$j(u) = \frac{u}{m}$.
Define $a \colon V \times V \to \Ci$ by
$a(u,v) = \int_\Omega \nabla u \, \overline{\nabla v}$.
Then $a$ is continuous and symmetric.
Since $\Omega$ is bounded it follows  from the 
(Dirichlet type) Poincar\'e inequality that
the norm 
\[
u \mapsto \Bigg( \int_\Omega |\nabla u|^2 + \int_\Omega \frac{|u|^2}{m^2} \Bigg)^{1/2}
\]
is an equivalent norm on $V$.
Therefore the form $a$ is $j$-elliptic.
Let $A$ be the operator associated with $(a,j)$.
We shall show that $A = -(m \Delta m)$.

Let $w \in D(A)$ and write $f = A w$.
Then there exists a $u \in V$ such that $w = j(u) = \frac{u}{m}$
and 
$\int_\Omega \nabla u \, \overline{\nabla v} = \int_\Omega f \, \frac{\overline{v}}{m}$
for all $v \in V$.
Observe that $\frac{f}{m} \in L_{1,\loc}(\Omega)$.
Taking $v \in \cd(\Omega)$ one deduces that 
$- \Delta u = \frac{f}{m}$ in $\cd(\Omega)'$.
Thus $w \in D((m \Delta m))$ and $- (m \Delta m) \, w = f$.

Conversely, let $w \in D((m \Delta m))$ and write $f = - (m \Delta m) \, w$.
Set $u = m \, w \in H^1_0(\Omega)$.
Then 
\[
a(u,v)
= \int_\Omega \nabla u \, \overline{\nabla v}
= - \langle \Delta u, \overline v \rangle
= \langle \frac{f}{m}, \overline v \rangle
= \int_\Omega f \, \frac{\overline{v}}{m} 
= \int_\Omega f \, \overline{j(v)}
\]
for all $v \in \cd(\Omega)$.
Since $\cd(\Omega)$ is dense in $V$ by 
\cite{AC}, Proposition~3.2,
it follows that $a(u,v) = \int_\Omega f \, \overline{j(v)}$ for all $v \in V$.
Thus $w = j(u) \in D(A)$.
This proves that $A = - (m \Delta m)$.

The operator $A$ is self-adjoint since $a$ is symmetric.
We next show the invariance of the set $C$.
The set $C$ is closed and convex in $L_2(\Omega)$.
Define $P \colon L_2(\Omega) \to C$ by 
$P f = (\RRe f) \wedge \frac{1}{m}$.
The $P$ is the orthogonal projection onto $C$.
Let $u \in V$.
Define $w = (\RRe u) \wedge 1 \in V$.
Then $P j(u) = j(w)$ and $\RRe a(w,u-w) = 0$.
Hence it follows from Proposition~\ref{psesqui202} that the set 
$C$ is invariant under $S$.
Since $f \leq 0$ if and only if $nf \in C$ for all $n \in \Ni$ the 
invariance of $C$ also implies that the semigroup is positive.\hfill$\Box$

\ruimte

By a similarity transformation we obtain two further kinds of 
multiplicative perturbations.
We leave the proofs to the reader.

\begin{prop} \label{psesqui466}
Let $\rho \colon \Omega \to (0,\infty)$ be such that 
$\frac{1}{\rho} \in L_{1,\loc}(\Omega)$.
Define the operator, formally denoted by $(\rho \Delta)$ on 
$L_2(\Omega, \frac{1}{\rho} \, dx)$ by the following.
Let $w,f \in L_2(\Omega, \frac{1}{\rho} \, dx)$.
Then we define $w \in D((\rho \Delta))$ and $(\rho \Delta) \, w = f$ if and only if
$w \in L_2(\Omega, \frac{1}{\rho} \, dx) \cap H^1_0(\Omega)$ 
and $\Delta w = \frac{f}{\rho}$ in $\cd(\Omega)'$.

Then the operator $(\rho \Delta)$ is self-adjoint and generates a 
submarkovian semigroup.
\end{prop}

\begin{prop} \label{psesqui467}
Let $\rho \colon \Omega \to (0,\infty)$ be such that 
$\frac{1}{\rho} \in L_{1,\loc}(\Omega)$.
Define the operator, formally denoted by $(\Delta \rho)$ on 
$L_2(\Omega, \rho \, dx)$ by the following.
Let $w,f \in L_2(\Omega, \rho \, dx)$.
Then $w \in D((\Delta \rho))$ and $(\Delta \rho) \, w = f$ if and only if
$\rho \, w \in H^1_0(\Omega)$ 
and $\Delta (\rho \, w) = f$ in $\cd(\Omega)'$.

Then the operator $(\Delta \rho)$ is self-adjoint and generates
a submarkovian semigroup.
\end{prop}

\subsection{Robin boundary conditions} \label{Ssesqui4.5}

Let $\Omega \subset \Ri^d$ be an open set with arbitrary 
boundary $\Gamma$.
At first we consider an arbitrary Borel measure on $\Gamma$ and then 
specialize to the $(d-1)$-dimensional Hausdorff measure.

For all $i,j \in \{ 1,\ldots,d \} $ let $a_{ij} \in L_\infty(\Omega,\Ci)$.
Let $\theta \in [0,\frac{\pi}{2})$.
Suppose 
$\sum_{i,j=1}^d a_{ij}(x) \, \xi_i \, \overline{\xi_j} \in \Sigma_\theta$
for all $\xi \in \Ci^d$ and a.e.\ $x \in \Omega$.
Let $\mu$ be a (positive) Borel measure on $\Gamma$ such that $\mu(K) < \infty$
for every compact $K \subset \Gamma$.
Define the form $a$ by 
\[
D(a) 
= \{ u \in H^1(\Omega) \cap C(\overline \Omega) 
     : \int_\Gamma |u|^2 \, d\mu < \infty \} 
\]
and 
\[
a(u,v) 
= \sum_{i,j=1}^d \int_\Omega (\partial_i u) \, a_{ij} \, \overline{\partial_j v}
   + \int_\Gamma u \, \overline v \, d\mu
.  \]
Then $C_c^\infty(\Omega) \subset D(a) \subset L_2(\Omega)$ and $a$ is sectorial.
In order to characterize the associated operator $A$ we
need to introduce two concepts and one more condition.
First, define the Neumann form $a_N$ by 
$D(a_N) = H^1(\Omega)$ and
\[
a_N(u,v) 
= \sum_{i,j=1}^d \int_\Omega (\partial_i u) \, a_{ij} \, \overline{\partial_j v}
.  \]
Throughout this subsection we suppose the form $a_N$ is closable.
Here we are more interested in the degeneracy caused by~$\mu$.
If $u \in D(\overline{a_N})$ and $f \in L_2(\Omega)$, then we say that 
$\ca u = f$ {\bf weakly on $\Omega$} if 
\[
\overline{a_N}(u,v) = \int_\Omega f \, \overline v
\]
for all $v \in C_c^\infty(\Omega)$.
If $u \in D(\overline{a_N})$, then we say that $\ca u \in L_2(\Omega)$ 
{\bf weakly on $\Omega$}
if there exists an $f \in L_2(\Omega)$ such that $\ca u = f$ weakly on $\Omega$.
Clearly such a function $f$ is unique, if it exists.
Secondly, if $u \in D(\overline{a_N})$ and $\varphi \in L_2(\Gamma,\mu)$,
then we say that $\varphi$ is {\bf an $(a,\mu)$-trace of $u$}, or shortly,
{\bf a trace of $u$}, if there exist $u_1,u_2,\ldots \in D(a)$
such that $\lim u_n = u$ in $D(\overline{a_N})$ and $\lim u_n|_\Gamma = \varphi$
in $L_2(\Gamma,\mu)$.
Moreover, let $H^1_{a,\mu}(\Omega)$ be the set of all $u \in D(\overline{a_N})$
for which there exists a $\varphi \in L_2(\Gamma,\mu)$ such that $\varphi$ is 
a trace of $u$.
We emphasize that $\varphi$ is not unique (almost everywhere) in general.
Clearly $D(a) \subset H^1_{a,\mu}(\Omega)$.
With the help of these definitions we can describe the operator~$A$ as follows.

\begin{prop} \label{psesqui450}
Let $u,f \in L_2(\Omega)$. 
Then $u \in D(A)$ and $A u = f$ if and only if
$u \in H^1_{a,\mu}(\Omega)$, 
$\ca u = f$ weakly on $\Omega$ and there exists a 
$\varphi \in L_2(\Gamma, \mu)$
such that $\varphi$ is a trace of $u$ and 
\begin{equation}
\overline{a_N}(u,v) - \int_\Omega (\ca u) \, \overline v 
= - \int_\Gamma \varphi \, \overline v \, d\mu
\label{epsesqui450;1}
\end{equation}
for all $v \in D(a)$.

If the conditions are valid, then the function $\varphi$ is unique.
\end{prop}
\proof\
`$\Rightarrow$'.
There exists a Cauchy sequence $u_1,u_2,\ldots$ in $D(a)$ 
such that $\lim u_n = u$ in $L_2(\Omega)$ and 
$\lim a(u_n,v) = (f,v)_H$ for all $v \in D(a)$.
Then $u_1,u_2,\ldots$ is a Cauchy sequence in $D(\overline{a_N})$.
Therefore $u \in D(\overline{a_N})$ and $\lim u_n = u$ in $D(\overline{a_N})$.
Moreover, $u_1|_\Gamma, u_2|_\Gamma, \ldots$ is a Cauchy sequence in 
$L_2(\Gamma,\mu)$.
Therefore $\varphi := \lim u_n|_\Gamma$ exists in $L_2(\Gamma,\mu)$.
Then $\varphi$ is a trace of $u$.
Let $v \in D(a)$.
Then 
\[
\overline{a_N}(u,v) + \int_\Gamma \, \varphi \, \overline v \, d\mu
= \lim a(u_n,v)
= (f,v)_H
= \int_\Omega f \, \overline v
.  \]
Therefore if  $v \in C_c^\infty(\Omega)$, then 
\[
\overline{a_N}(u,v)
= \int_\Omega f \, \overline v
,  \]
so $\ca u = f$ weakly on $\Omega$.
Moreover, 
\[
\overline{a_N}(u,v) + \int_\Gamma \, \varphi \, \overline v \, d\mu
= \int_\Omega (\ca u) \, \overline v
\]
for all $v \in D(a)$.

If also $\varphi' \in L_2(\Gamma,\mu)$ satisfies (\ref{epsesqui450;1}),
then $\int_\Gamma (\varphi - \varphi') \, \overline v \, d\mu = 0$ for 
all $v \in D(a)$. 
But the space $ \{ v|_\Gamma : v \in H^1(\Omega) \cap C_c(\overline \Omega) \} $
is a $*$-algebra which separates the points of $\Gamma$.
Therefore it is dense in $C_0(\Gamma)$.
Let $\psi \in C_c(\Gamma)$.
Then there exists a $\chi \in C_c^\infty(\Ri^d)$ such that 
$\chi|_{\supp \psi} = \one$.
By the above there exist $v_1,v_2,\ldots \in H^1(\Omega) \cap C_c(\overline \Omega)$
such that $\lim v_n|_\Gamma = \psi$ in $C_0(\Gamma)$.
Then $\lim (\chi \, v_n)|_\Gamma = \psi$ in $C_0(\Gamma)$.
Moreover, $\mu(\supp (\chi|_\Gamma)) < \infty$.
Therefore $\lim (\chi \, v_n)|_\Gamma = \psi$ in $L_2(\Gamma,\mu)$
and the space $ \{ v|_\Gamma : v \in H^1(\Omega) \cap C_c(\overline \Omega) \} $
is dense in $L_2(\Gamma,\mu)$.
Thus $\varphi' = \varphi$.

`$\Leftarrow$'.
There exist $\varphi \in L_2(\Gamma,\mu)$ and 
$u_1,u_2,\ldots \in D(a)$ such that 
$\lim u_n = u$ in $D(\overline{a_N})$,
$\lim u_n|_\Gamma = \varphi$ in $L_2(\Gamma,\mu)$
and (\ref{epsesqui450;1}) is valid for all 
$v \in D(a)$.
Then $u_1,u_2,\ldots$ is a Cauchy sequence in $D(a)$
and 
\[
\lim_{n \to \infty} a(u_n,v)
= \overline{a_N}(u,v)
   + \int_\Gamma \, \varphi \, \overline v \, d\mu
= \int_\Omega (\ca u) \, \overline v
= \int_\Omega f \, \overline v
\]
for all $v \in D(a)$.
So $u \in D(A)$ and $A u = f$.\hfill$\Box$

\ruimte

This proposition shows how our general results can be easily applied.
It is worthwhile to consider more closely the associated closed form since
this is intimately related to the problem to define a trace
in $L_2(\Gamma,\mu)$ of suitable functions in $H^1(\Omega)$.

Let
\[
W = \overline{ \{ (u, u|_\Gamma) : u \in D(a) \} }
,  \]
where the closure is in $D(\overline{a_N}) \oplus L_2(\Gamma,\mu)$.
Then the map $u \mapsto (u, u|_\Gamma)$ from $D(a)$ into $W$
is an isometry and therefore it extends to a unitary map
from the completion of $D(a)$ onto~$W$.
The form $a$ is closable if and only if the map $j \colon W \to L_2(\Omega)$ 
defined by $j(u,\varphi) = u$
is injective.
Note that if $\varphi \in L_2(\Gamma,\mu)$, then
$(0,\varphi) \in W$ if and only if $\varphi$ is a trace of~$0$.

The following lemma is due to Daners \cite{Daners2} Proposition~3.3
in the strongly elliptic case, but our proof is different.

\begin{lemma} \label{lsesqui450}
There exists a Borel set $\Gamma_{a,\mu} \subset \Gamma$ such that 
\[
\{ \varphi \in L_2(\Gamma,\mu) : \varphi \mbox{ is a trace of } 0 \} 
= L_2(\Gamma \setminus \Gamma_{a,\mu}, \mu)
.  \]
\end{lemma}
\proof\
Set $F = \{ \varphi \in L_2(\Gamma,\mu) : (0,\varphi) \in W \} $.
Then $F$ is a closed subspace of $L_2(\Gamma,\mu)$.

First we show that $u \, \psi \in F$ for all $\psi \in F$ and 
$u \in D(a) \cap W_\infty^1(\Ri^d)$.
Since $\psi \in F$ there exist $u_1,u_2,\ldots \in D(a)$ 
such that $\lim u_n = 0$ in $D(\overline{a_N})$ and 
$\lim u_n|_\Gamma = \psi$ in $L_2(\Gamma,\mu)$.
Then $u \, u_n \in D(a)$ for all $n \in \Ni$ and
$\lim (u \, u_n)|_\Gamma = u \, \psi$ in $L_2(\Gamma,\mu)$.
By the Leibniz rule one deduces that 
\[
(\Re \overline{a_N})(u \, u_n)^{1/2} 
\leq 
   \|u_n\|_2 \, \Big( \sum \Big\| \, |\frac{a_{ij} + \overline{a_{ji}}}{2}| 
                        \, |\partial_i u| \, |\partial_j u| \, \Big\|_\infty \Big)^{1/2}
  + \|u\|_\infty \, (\Re \overline{a_N})(u_n)^{1/2} 
\]
for all $n \in \Ni$ and $\lim u \, u_n = 0$ in $D(\overline{a_N})$.
So $u \, \psi \in F$.

Secondly, let $P \colon L_2(\Gamma,\mu) \to F$ be the orthogonal projection.
Let $\varphi \in L_2(\Gamma,\mu)$ and suppose that $\mu([\varphi \neq 0]) < \infty$.
We shall prove that $P \varphi = 0$ a.e.\ on $[\varphi = 0]$.
Let $A = [\varphi \neq 0]$.
Since 
$ \{ u|_\Gamma : u \in H^1(\Omega) \cap C_c^\infty(\Ri^d) \} $ is dense
in $L_2(\Gamma,\mu)$ there exist 
$u_1,u_2,\ldots \in H^1(\Omega) \cap C_c^\infty(\Ri^d)$
such that $\lim u_n|_\Gamma = \one_A$ in $L_2(\Gamma,\mu)$.
Then also $\lim (0 \vee \RRe u_n \wedge \one)|_\Gamma = \one_A$ in $L_2(\Gamma,\mu)$,
so we may assume that $u_n \in D(a) \cap W_\infty^1(\Ri^d)$ and 
$0 \leq u_n \leq \one$ for all $n \in \Ni$.
Passing to a subsequence if necessary we may assume that 
$\lim u_n|_\Gamma = \one_A$ a.e.
Therefore $\lim u_n \, P \varphi = \one_A \, P \varphi$ in $L_2(\Gamma,\mu)$.
Since $u_n \, P \varphi \in F$ for all $n \in \Ni$ one deduces that 
$\one_A \, P \varphi \in F$.
Then 
\[
\|\varphi - \one_A \, P \varphi\|_{L_2(\Gamma,\mu)}
= \|\one_A( \varphi - P \varphi)\|_{L_2(\Gamma,\mu)}
\leq \|\varphi - P \varphi\|_{L_2(\Gamma,\mu)}
.  \]
So $\one_A \, P \varphi = P \varphi$ and 
$P \varphi = 0$ a.e.\ on $A^{\rm c} = [\varphi = 0]$.
Now the lemma easily follows from Zaanen's theorem 
\cite{ArT} Proposition~1.7.\hfill$\Box$

\ruimte

Obviously the set $\Gamma_{a,\mu}$ in Lemma~\ref{lsesqui450} is unique in 
the sense that $\mu(\Gamma_{a,\mu} \Delta \Gamma') = 0$ whenever 
$\Gamma' \subset \Gamma$ is another Borel set 
with this property.
It follows from the last paragraph of Section~3 in \cite{AW2} that the 
set $\Gamma_{a,\mu}$ coincides with the set $S$ in \cite{AW2} Proposition~3.6.
In \cite{AW2} Example~4.2 there is an example of an open set $\Omega$ 
and $\mu$ the $(d-1)$-dimensional Hausdorff measure such that 
$\mu(\Gamma) < \infty$ and $\mu(\Gamma \setminus \Gamma_{a,\mu}) > 0$,
where $a$ is the form of the Laplacian.

It is clear from the construction of 
$\Gamma_{a,\mu}$ and definition of $H^1_{a,\mu}(\Omega)$ that 
there exists a unique map $\Tr_{a,\mu} \colon H^1_{a,\mu}(\Omega) \to L_2(\Gamma_{a,\mu},\mu)$
in a natural way, which we call {\bf trace}.
Note that if $u \in H^1_{a,\mu}(\Omega)$, then $\Tr_{a,\mu} u$ 
is the unique $\varphi \in L_2(\Gamma_{a,\mu},\mu)$ such that $\varphi$ is 
an $(a,\mu)$-trace of $u$.
In general, however, the map $\Tr_{a,\mu}$ is not continuous from 
$(H^1_{a,\mu}(\Omega),\|\cdot\|_{\overline{a_N}})$ into 
$L_2(\Gamma_{a,\mu},\mu)$.
A counter-example is in \cite{Daners2}, Remark~3.5(f).

The map $u \mapsto (u,\Tr_{a,\mu} u)$ from $H^1_{a,\mu}(\Omega)$ into 
$D(\overline{a_N}) \oplus L_2(\Gamma_{a,\mu},\mu)$ is injective.
Therefore one can define a norm on 
$H^1_{a,\mu}(\Omega)$ by
\[
\|u\|_{H^1_{a,\mu}(\Omega)}^2
= \|u\|_{D(\overline{a_N})}^2 + \|\Tr_{a,\mu} u\|_{L_2(\Gamma_{a,\mu},\mu)}^2
.  \]
It is easy to verify that $H^1_{a,\mu}(\Omega)$ is a Hilbert space.
Moreover, the map $\Tr_{a,\mu} \colon H^1_{a,\mu}(\Omega) \to L_2(\Gamma_{a,\mu},\mu)$ is a continuous 
linear operator with dense range.

It is now possible to reconsider the element $\varphi \in L_2(\Gamma,\mu)$
in Proposition~\ref{psesqui450}.

\begin{prop} \label{psesqui459}
Let $u,f \in L_2(\Omega)$. 
Then $u \in D(A)$ and $A u = f$ if and only if
$u \in H^1_{a,\mu}(\Omega)$, 
$\ca u = f$ weakly on $\Omega$ and 
\[
\overline{a_N}(u,v) - \int_\Omega (\ca u) \, \overline v 
= - \int_\Gamma (\Tr_{a,\mu} u) \, \overline v \, d\mu
\]
for all $v \in D(a)$.
\end{prop}
\proof\
Let $u \in D(A)$ and $\varphi \in L_2(\Gamma,\mu)$ be the 
corresponding unique
element as in Proposition~\ref{psesqui450}.
If $\psi \in L_2(\Gamma \setminus \Gamma_{a,\mu},\mu) = F$, then 
there exist $v_1,v_2,\ldots \in D(a)$
such that $\lim v_n = 0$ in $D(\overline{a_N})$ and 
$\lim v_n|_\Gamma = \psi$ in $L_2(\Gamma,\mu)$.
Substituting $v = v_n$ in (\ref{epsesqui450;1}) and taking the 
limit $n \to \infty$ one deduces that 
$\int_\Gamma \varphi \, \overline{\psi} \, d\mu = 0$.
So $\varphi \in L_2(\Gamma_{a,\mu},\mu)$ and $\varphi = \Tr_{a,\mu} u$.\hfill$\Box$

\ruimte

We now consider the case where $\mu$ is the $(d-1)$-dimensional Hausdorff measure,
which we denote by $\sigma$.
In particular, we assume that $\sigma(K) < \infty$ for every compact $K \subset \Gamma$.
Moreover, we write $\Gamma_a = \Gamma_{a,\sigma}$ and $\Tr_a = \Tr_{a,\sigma}$.
The measure $\sigma$ coincides with the usual surface measure if $\Omega$ is $C^1$.
We continue to consider, however, the case where $\Omega$ is an arbitrary
bounded open set.
If $\Omega$ has a Lipschitz 
continuous boundary and the form $a$ equals the 
the classical Dirichlet form $l$, then $\Gamma_l = \Gamma_a = \Omega$ 
by the trace theorem (see \cite{Nec}
Th\'eor\`eme~2.4.2).
By \cite{AW2} Proposition~5.5 it follows that $\sigma(\Gamma_l) > 0$ if $\Omega$ is 
bounded, without any regularity condition on the boundary.
(Note, however, that there exists an open connected subset $\Omega \subset \Ri^3$
such that $\sigma(\Gamma \setminus \Gamma_l) > 0$,
see \cite{AW2}, Example~4.3).
The embedding of $H^1_{l,\sigma}(\Omega)$ into $L_2(\Omega)$ is compact 
if $\Omega$ has finite measure, by \cite{AW2} Corollary~5.2.
This surprising phenomenon is a consequence
of Maz'ya's inequality.
It was Daners \cite{Daners2} who was the first to exploit this inequality to establish
results for Robin boundary conditions on rough domains.
Further results were given in \cite{AW2} Section~5.

We conclude our remarks by considering $\mu = \beta \, \sigma$,
where $\beta \in L_\infty(\Gamma,\Ri)$ and $\beta \geq 0$ a.e.
We define the  weak normal derivative with respect to 
the matrix $(a_{ij})$.
Let $\varphi \in L_2(\Gamma,\mu)$, $u \in D(\overline{a_N})$ and suppose that 
$\ca u \in L_2(\Omega)$ weakly on $\Omega$.
Then we say that $\varphi$ is the {\bf $(a_{ij})$-normal derivative of $u$}
if 
\[
\overline{a_N}(u,v) - \int_\Omega (\ca u) \, \overline v 
= \int_\Gamma \varphi \, \overline v \, d\sigma
\]
for all $v \in D(a)$.
If $\Omega$ is of class $C^1$, $\mu$ is the $(d-1)$-dimensional Hausdorff measure
and $u \in C^1(\overline \Omega)$, then 
our weak definition coincides with the classical definition
by Green's theorem.
We reformulate Proposition~\ref{psesqui450}.

\begin{prop} \label{psesqui452}
Let $u,f \in L_2(\Omega)$. 
Then $u \in D(A)$ and $A u = f$ if and only if
$u \in H^1_{a,\beta \, \sigma}(\Omega)$, 
$\ca u = f$ weakly on $\Omega$ and 
$- \beta \, \Tr_{a,\beta \, \sigma} u$ is the $(a_{ij})$-normal derivative of $u$.
\end{prop}

Note that if  $(a_{ij})$ 
is strongly elliptic and if $u \in D(A)$ and $A u = f$,
then $u \in H^1(\Omega)$, $\ca u = f$ weakly on $\Omega$, $u$ has a trace $\Tr u$
and $\nu \cdot a \nabla u = - \beta \Tr u$ weakly.
Thus one recovers the classical statement.

\subsection{The Dirichlet-to-Neumann operator} \label{Ssesqui4.2}

Let $\Omega$ be a bounded open subset of $\Ri^d$ with Lipschitz boundary~$\Gamma$,
provided with the $(d-1)$-dimensional Hausdorff measure.
Let $\Tr \colon H^1(\Omega) \to L_2(\Gamma)$ be the trace map.
We denote by $\Delta_D$ the Dirichlet Laplacian on $\Omega$.
If $\varphi \in L_2(\Gamma)$, $u \in H^1(\Omega)$ and 
$\Delta u \in L_2(\Omega)$ weakly on $\Omega$, then we say that 
$\frac{\partial u}{\partial \nu} = \varphi$ weakly if 
$\varphi$ is the $(a_{ij})$-normal derivative of $u$, where 
$a_{ij} = \delta_{ij}$.

Let $\lambda \in \Ri$ and suppose that $\lambda \not\in \sigma(-\Delta_D)$.
The {\bf Dirichlet-to-Neumann operator} $D_\lambda$ on $L_2(\Gamma)$ is defined
as follows.
Let $\varphi,\psi \in L_2(\Gamma)$. 
Then we define $\varphi \in D(D_\lambda)$ and $D_\lambda \varphi = \psi$
if there exists a $u \in H^1(\Omega)$ such that 
$-\Delta u = \lambda u$ weakly on $\Omega$, $\Tr u = \varphi$ and 
$\frac{\partial u}{\partial \nu} = \psi$ weakly.
We next show that the Dirichlet-to-Neumann operator is an example
of the $m$-sectorial operators obtained in Corollary~\ref{csesqui201.2}.

Define the sesquilinear form $a \colon H^1(\Omega) \times H^1(\Omega) \to \Ci$ 
by 
\[
a(u,v) 
= \int_\Omega \nabla u \, \overline {\nabla v} 
   - \lambda \int_\Omega u \, \overline v
.  \]
Moreover, define $j \colon H^1(\Omega) \to L_2(\Gamma)$ by 
$j(u) = \Tr u$.
Then $\ker j = H^1_0(\Omega)$ by \cite{Alt} Lemma~A.6.10.
Clearly the form $a$ is continuous, the map $j$ is bounded and $j(H^1(\Omega))$ is 
dense in $L_2(\Gamma)$.
Using the definitions one deduces that 
\[
V(a)
= \{ u \in H^1(\Omega) : -\Delta u = \lambda u \mbox{ weakly on } \Omega \} 
.  \]
It follows from Step~1 in the proof of Proposition~3.3 in \cite{ArM}
that there exist $\omega \in \Ri$ and $\mu > 0$ such that 
\[
\RRe a(u) + \omega \, \|j(u)\|_{L_2(\Gamma)}^2 \geq \mu \, \|u\|_V^2
\]
for all $u \in V(a)$.
Moreover, $H^1(\Omega) = V(a) + \ker j$ by Lemma~3.2 in \cite{ArM}.
So the conditions of Corollary~\ref{csesqui201.2} are satisfied.

Note that if $\lambda_1$ is the lowest eigenvalue of the operator
$-\Delta_D$ on $\Omega$ and $u \in H^1_0(\Omega)$ is an eigenfunction with 
eigenvalue $\lambda_1$, then 
(\ref{etsesqui201;1}) is not valid if $\lambda > \lambda_1$. 
Therefore Theorem~\ref{tsesqui201} is not applicable and this 
example is the reason why we used the space $V(a)$ in Corollary~\ref{csesqui201.2}.
If still $\lambda > \lambda_1$ then there exists a $\tilde u \in H^1_0(\Omega)$
such that $\tilde u \neq 0$ and $a(\tilde u) = 0$.
If $U = V(a) + \spann \{ \tilde u \} $ then $U$ is a closed subspace of 
$H^1(\Omega)$ such that $D_H(a) \subset U$.
But $(a|_{U \times U}, j|_U)$ does not satisfy (\ref{ecsesqui201.2;1}).
Hence Proposition~\ref{psesqui204}\ref{psesqui204-2} cannot be extended 
to the setting of Corollary~\ref{csesqui201.2}.

Let $A$ be the operator associated with $(a,j)$.
We next show that $A = D_\lambda$.
Let $\varphi,\psi \in L_2(\Gamma)$. 
Suppose $\varphi \in D(A)$ and $A \varphi = \psi$.
Then there is a $u \in H^1(\Omega)$ such that $\Tr u = \varphi$
and $a(u,v) = (\psi, \Tr v)_{L_2(\Gamma)}$ for all $v \in H^1(\Omega)$.
For all $v \in H^1_0(\Omega)$ one has 
\[
\int_\Omega \nabla u \, \overline {\nabla v} 
   - \lambda \int_\Omega u \, \overline v
= a(u,v) = 0
,  \]
so $-\Delta u = \lambda u$ weakly on $\Omega$.
Then 
\[
\int_\Omega \nabla u \, \overline {\nabla v} 
   + \int_\Omega (\Delta u) \, \overline v
= a(u,v)
= (\psi, \Tr v)_{L_2(\Gamma)}
\]
for all $v \in H^1(\Omega)$.
So $\frac{\partial u}{\partial \nu} = \psi$ weakly.
Therefore $\varphi \in D(D_\lambda)$ and $D_\lambda \varphi = \psi$.
Conversely, suppose $\varphi \in D(D_\lambda)$ and $D_\lambda \varphi = \psi$.
By definition there exists a $u \in H^1(\Omega)$ such that 
$-\Delta u = \lambda u$ weakly on $\Omega$, $\Tr u = \varphi$ and 
$\frac{\partial u}{\partial \nu} = \psi$ weakly.
Then 
\[
a(u,v)
= \int_\Omega \nabla u \, \overline {\nabla v} 
   - \lambda \int_\Omega u \, \overline v
= \int_\Omega \nabla u \, \overline {\nabla v} 
   + \int_\Omega (\Delta u) \, \overline v
= \int_\Gamma \frac{\partial u}{\partial \nu} \, \overline{\Tr v}
= (\psi, \Tr v)_{L_2(\Gamma)}
\]
for all $v \in H^1(\Omega)$.
So $\varphi = j(u) \in D(A)$ and $\psi = A j(u) = A \varphi$.
Thus $D_\lambda = A$ is the operator associated with $(a,j)$.

If $S$ is the semigroup generated by $- D_\lambda$,
then it follows as in the proof of Corollary~\ref{csesqui331} that 
$S$ is real and positive.
Moreover, if $\lambda \leq 0$, then $S$ extends consistently to a 
continuous contraction semigroup
on $L_p(\Omega)$ for all $p \in [1,\infty]$.

\subsection{Wentzell boundary conditions} \label{Ssesqui4.3}

Let again $\Omega$ be an open subset of $\Ri^d$ with arbitrary boundary~$\Gamma$
and let $\sigma$ be the $(d-1)$-dimensional Hausdorff measure on $\Gamma$.
We assume that $\sigma(K) < \infty$ for every compact $K \subset \Gamma$.
All $L_p$ spaces on $\Gamma$ are with respect to the measure
$\sigma$, except if written different explicitly.
For all $i,j \in \{ 1,\ldots,d \} $ let $a_{ij} \in L_\infty(\Omega)$.
Let $\theta \in [0,\frac{\pi}{2})$.
Suppose 
$\sum_{i,j=1}^d a_{ij}(x) \, \xi_i \, \overline{\xi_j} \in \Sigma_\theta$
for all $\xi \in \Ci^d$ and a.e.\ $x \in \Omega$.
Define the form $b$ by 
\[
D(b) 
= \{ u \in H^1(\Omega) \cap C(\overline \Omega) 
     : \int_\Gamma |u|^2 \, d\sigma < \infty \} 
\]
and
\[
b(u,v) 
= \sum_{i,j=1}^d \int_\Omega (\partial_i u) \, a_{ij} \, \overline{\partial_j v}
   + \int_\Gamma  \, u \, \overline v \, d\sigma
.  \]
As in Subsection~\ref{Ssesqui4.5} we define the Neumann form $b_N$ by 
$D(b_N) = H^1(\Omega)$ and
\[
b_N(u,v) 
= \sum_{i,j=1}^d \int_\Omega (\partial_i u) \, a_{ij} \, \overline{\partial_j v}
.  \]
Throughout this subsection we assume that the form $b_N$ is closable.
Set $\widetilde \Gamma = \Gamma_{b,\sigma}$ and $\Tr = \Tr_{b,\sigma}$.
Moreover, we assume that the map 
$\Tr \colon (H^1_{b,\sigma}(\Omega), \|\cdot\|_{\overline{b_N}}) \to L_2(\widetilde \Gamma)$ is 
continuous.

Fix $\alpha \in L_\infty(\widetilde \Gamma)$ and 
$B \in \cl(L_2(\widetilde \Gamma))$.
Throughout this subsection we assume that there exists an $\omega > 0$ such that 
\begin{equation}
\omega \, \|B \varphi\|_{L_2(\widetilde \Gamma)}^2
+ \int_{\widetilde \Gamma} \RRe \alpha \, |\varphi|^2
\geq 0
\label{eSsesqui4;21}
\end{equation}
for all $\varphi \in L_2(\widetilde \Gamma)$.
As an example, if $\beta \in L_\infty(\Gamma)$ and $B$ is the 
multiplication operator with~$\beta$, then the assumption
\[
\omega \, |\beta|^2 + \RRe \alpha \geq 0
\]
for some $\omega > 0$, implies (\ref{eSsesqui4;21}).

Define the form $a$ by 
\[
D(a) = H^1_{b,\sigma}(\Omega)
\]
and 
\[
a(u,v) 
= \sum_{i,j=1}^d \int_\Omega (\partial_i u) \, a_{ij} \, \overline{\partial_j v}
   + \int_{\widetilde \Gamma} \Tr u \, \overline{\Tr v} \, \alpha \, d\sigma
.  \]
Then $a$ is continuous.
Let $H$ be the closure of the space
$ \{ (u,B(\Tr u)) : u \in H^1_{b,\sigma}(\Omega) \} $ in 
the space $L_2(\Omega) \oplus L_2(\widetilde \Gamma)$ with 
induced norm.
Define the injective map $j \colon H^1_{b,\sigma}(\Omega) \to H$ by 
\[
j(u) = (u, B(\Tr u))
.  \]
If $B$ has dense range, then $H = L_2(\Omega) \oplus L_2(\widetilde \Gamma)$
since the space $  \{ (u,\Tr u) : u \in H^1_{b,\sigma}(\Omega) \}$
is dense in $L_2(\Omega) \oplus L_2(\widetilde \Gamma)$ by Step a) in the 
proof of Theorem~2.3 in \cite{AMPR}.
Then the claim follows by the range condition on $B$.
Note that the condition~(\ref{eSsesqui4;21}) 
together with the assumed continuity of 
$\Tr \colon (H^1_{b,\sigma}(\Omega), \|\cdot\|_{\overline{b_N}}) \to L_2(\widetilde \Gamma)$
imply that $a$ is $j$-elliptic.
Let $A$ be the operator associated with $(a,j)$.

\begin{prop} \label{psesqui470}
Let $x,y \in H$.
Then $x \in D(A)$ and $A x = y$ if and only if there exist
$u \in H^1_{b,\sigma}(\Omega)$ and $\psi \in L_2(\widetilde \Gamma)$ 
such that $x = (u, B(\Tr u))$,
$\ca u \in L_2(\Omega)$ weakly on $\Omega$,  
$(B^* \psi - \alpha \, \Tr u)$ is the $(a_{ij})$-normal derivative of $u$
and $y = (\ca u, \psi)$.
\end{prop}
\proof\
`$\Rightarrow$'.
Write $y = (f,\psi) \in H$.
There exists a $u \in H^1_{b,\sigma}(\Omega)$ such that $x = j(u)$ and 
\[
\overline{b_N}(u,v) 
   + \int_{\widetilde \Gamma} \Tr u \, \overline{\Tr v} \, \alpha \, d\sigma
= (y,j(v))_H
= \int_\Omega f \, \overline{v} + \int_{\widetilde \Gamma} \psi \, \overline{B(\Tr v)} \, d\sigma
\]
for all $v \in H^1_{b,\sigma}(\Omega)$.
Taking only $v \in C_c^\infty(\Omega)$ one deduces that 
$\ca u = f$ weakly on $\Omega$.
In particular, $y = (f,\psi) = (\ca u, \psi)$.
Moreover, 
\[
\overline{b_N}(u,v)
- \int_\Omega (\ca u) \, \overline v
= \int_{\widetilde \Gamma} (B^* \psi - \alpha \, \Tr u) \, \overline{\Tr v} \, d\sigma
\]
for all $v \in H^1_{b,\sigma}(\Omega)$, which implies that 
$(B^* \psi - \alpha \, \Tr u)$ is the $(a_{ij})$-normal derivative of $u$.

`$\Leftarrow$'.
Let $u \in H^1_{b,\sigma}(\Omega)$ and $\psi \in L_2(\widetilde \Gamma)$ 
be such that  $x = (u, B(\Tr u))$,
$\ca u \in L_2(\Omega)$ weakly on $\Omega$,  
$(B^* \psi - \alpha \, \Tr u)$ is the $(a_{ij})$-normal derivative of $u$
and $y = (\ca u, \psi)$.
Then $x = j(u)$.
Since $(B^* \psi - \alpha \, \Tr u)$ is the $(a_{ij})$-normal derivative of $u$
one deduces that  
\[
\overline{b_N}(u,v) - \int_{\widetilde \Gamma} (\ca u) \, \overline v 
= \int_{\widetilde \Gamma} (B^* \psi - \alpha \, \Tr u) \, \overline{\Tr v} \, d\sigma
\]
for all $v \in H^1_{b,\sigma}(\Omega)$.
So
\begin{eqnarray*}
a(u,v)
& = & \overline{b_N}(u,v) +  \int_{\widetilde \Gamma} \Tr u \, \overline{\Tr v} \, \alpha \, d\sigma  \\
& = & \int_\Omega (\ca u) \, \overline v 
   + \int_{\widetilde \Gamma} \psi \, \overline{B(\Tr v)} \, d\sigma
= (y,j(v))_H
\end{eqnarray*}
for all $v \in H^1_{b,\sigma}(\Omega)$.
Therefore $x \in D(A)$ and $A x = y$.\hfill$\Box$

\ruimte

Suppose that $B$ has dense range.
Then $H$ is isomorphic with $L_2(\Omega \sqcup \widetilde \Gamma)$ in a natural way,
where $\sqcup$ denotes the disjoint union of the measure spaces.
We use this isomorphism to identify $H$ with $L_2(\Omega \sqcup \widetilde \Gamma)$.
It is easy to verify as in the proof of 
Corollary~\ref{csesqui331}\ref{csesqui331-1}
that $S$ leaves $L_2(\Omega,\Ri) \oplus L_2(\widetilde \Gamma,\Ri)$ invariant
if the form $b_N$ is real, $\alpha$ is real valued and 
$B$ maps $L_2(\widetilde \Gamma ;\Ri)$ into itself.
We next characterize positivity of~$S$.

\begin{prop} \label{psesqui471}
Suppose the form $b_N$ is real, $\alpha$ is real valued
and $B$ maps $L_2(\widetilde \Gamma;\Ri)$ densely into itself.
\begin{tabel}
\item \label{psesqui471-1}
The map $B$ is a lattice homomorphism if and only if the semigroup $S$ is positive.
\item \label{psesqui471-2}
If $\sigma(\widetilde \Gamma) < \infty$, the map $B$ is a lattice homomorphism,
$\alpha \geq 0$ and there 
exists a $c \geq 1$ such that 
$\frac{1}{c} \, \one \leq B \one \leq c \, \one$, 
then $S$ 
extends continuously to a bounded semigroup on $L_\infty(\Omega \sqcup \widetilde \Gamma)$.
\end{tabel}
\end{prop}
\proof\
`\ref{psesqui471-1}'.
Let $C = \{ (u,\varphi) \in H : u \geq 0 \mbox{ and } \varphi \geq 0 \} $.
Then $C$ is closed and convex in~$H$. 
Define $P \colon H \to C$ by $P(u,\varphi) = ((\RRe u)^+, (\RRe \varphi)^+)$.
Then $P$ is the orthogonal projection onto $C$.

`$\Rightarrow$'.
Let $u \in H^1_{b,\sigma}(\Omega)$.
Then $(\RRe u)^+ \in H^1_{b,\sigma}(\Omega)$ and 
\[
j((\RRe u)^+) 
= ((\RRe u)^+, B(\Tr((\RRe u)^+)))
= ((\RRe u)^+, (\RRe B(\Tr u))^+)
= P j(u)
\]
since $B$ is a lattice homomorphism.
Moreover, 
\[
\RRe a((\RRe u)^+, u - (\RRe u)^+) 
= a((\RRe u)^+, - (\RRe u)^-) 
= 0
.  \]
So $C$ is invariant under $S$ by Proposition~\ref{psesqui202}.

`$\Leftarrow$'.
If $S$ is positive, then $C$ is invariant under $S$.
Let $u \in H^1_{b,\sigma}(\Omega)$.
It follows from Proposition~\ref{psesqui202} that 
there exists a $w \in H^1_{b,\sigma}(\Omega)$ 
such that $P j(u) = j(w)$.
Then $((\RRe u)^+, (\RRe B(\Tr u))^+) = P j(u) = j(w) = (w, B(\Tr w))$.
Therefore $w = (\RRe u)^+$ and 
\[
(\RRe B(\Tr u))^+ 
= B(\Tr w) 
= B(\Tr ((\RRe u)^+))
= B((\RRe \Tr u)^+)
.  \]
This is for all $u \in H^1_{b,\sigma}(\Omega)$. 
Since $\Tr H^1_{b,\sigma}(\Omega)$ is dense in $L_2(\widetilde \Gamma)$ one deduces that 
$(B \varphi)^+ = B(\varphi^+)$ for all $\varphi \in L_2(\widetilde \Gamma,\Ri)$.
So $B$ is a lattice homomorphism.

`\ref{psesqui471-2}'.
Let $C = \{ (u,\varphi) \in H : u \leq \one \mbox{ and } \varphi \leq B \one \} $.
Then $C$ is closed and convex.
Define $P \colon H \to C$ by 
$P(u,\varphi) = ((\RRe u) \wedge \one, (\RRe \varphi) \wedge B \one)$.
Then $P$ is the orthogonal projection of $H$ onto $C$.
Let $u \in H^1_{b,\sigma}(\Omega)$. 
Define $w = (\RRe u) \wedge \one$. 
Then $w \in H^1_{b,\sigma}(\Omega)$ and 
$P j(u) = ((\RRe u) \wedge \one, (\RRe B(\Tr u)) \wedge B \one)
= ((\RRe u) \wedge \one, B (\Tr ((\RRe u) \wedge \one)))
= j(w)$.
Moreover, 
\begin{eqnarray*}
\RRe a(w,u-w)
& = & \RRe a((\RRe u) \wedge \one, i \IIm u + (\RRe u - \one)^+)
= a((\RRe u) \wedge \one, (\RRe u - \one)^+)  \\
& = & \int_{\widetilde \Gamma} \alpha \, \Tr ((\RRe u) \wedge \one) \, \Tr ((\RRe u - \one)^+)
= \int_{\widetilde \Gamma} \alpha  \, \Tr ((\RRe u - \one)^+)
\geq 0
\end{eqnarray*}
So by Proposition~\ref{psesqui202} the set $C$ is invariant under $S$.

Finally, let $(u,\varphi) \in H$ and suppose that $u \leq \one$ and $\varphi \leq \one$.
Then $\frac{1}{c} \, \varphi \leq B \one$ 
and $\frac{1}{c} \, (u,\varphi) \in C$.
Let $t > 0$ and write $(v,\psi) = S_t(u,\varphi)$.
Then $\frac{1}{c} \, (v,\psi) \in C$.
Hence $v \leq c \, \one$ and $\psi \leq c \, B \one \leq c^2 \, \one$.
So $S$ extends to a continuous semigroup on $L_\infty$ and 
$\|S_t\|_{\infty \to \infty} \leq c^2$ for all $t > 0$.\hfill$\Box$

\ruimte

Using the operator $A$ one can define another semigroup generator
which looks different.
If $u \in D(\overline{b_N})$, then we say that $\ca u \in H^1_{b,\sigma}(\Omega)$ 
{\bf weakly on $\Omega$}
if there exists an $f \in H^1_{b,\sigma}(\Omega)$ such that $\ca u = f$ weakly on $\Omega$.

The Laplacian with Wentzell boundary conditions can be realized 
in the Sobolev space~$H^1$.
This has been carried out in \cite{FGGOR} Theorem~2.1.
We generalize this approach for the elliptic operator~$A$.

\begin{prop} \label{psesqui472}
Define the operator $A_1$ on $H^1_{b,\sigma}(\Omega)$ 
by taking as domain $D(A_1)$ the set of all 
$u \in H^1_{b,\sigma}(\Omega)$ such that 
$\ca u \in H^1_{b,\sigma}(\Omega)$ weakly on $\Omega$
and $(B^* B (\Tr \ca u) - \alpha \, \Tr u)$ is the $(a_{ij})$-normal derivative
of $u$
and letting $A_1 u = \ca u$ for all $u \in D(A_1)$.
Then $- A_1$ generates a holomorphic semigroup on $H^1_{b,\sigma}(\Omega)$.
\end{prop}
\proof\
Let $a_c$ be the classical form associated with $(a,j)$
(see Theorem~\ref{tsesqui212}).
Then $A$ is associated with the closed sectorial form $a_c$.
Define the operator $A_0$ in $H$ by 
$D(A_0) = \{ w \in D(A) : A w \in D(a_c) \} $ and 
$A_0 w = A w$ for all $w \in D(A_0)$.
Then $-A_0$ generates a holomorphic semigroup in the Hilbert space
$(D(a_c), \|\cdot\|_{a_c})$.
The map $j \colon H^1_{b,\sigma}(\Omega) \to D(a_c)$ is a isomorphism of normed spaces.
Hence the operator $- j^{-1} A_0 j$ generates a holomorphic semigroup 
on $H^1_{b,\sigma}(\Omega)$.
Therefore it suffices to show that $A_1 = j^{-1} A_0 j$.

Let $u \in D(j^{-1} A_0 j)$.
Then $j(u) \in D(A)$, $A j(u) \in j(H^1_{b,\sigma}(\Omega))$ and $A_0 j(u) = A j(u)$.
It follows from Proposition~\ref{psesqui470} that $\ca u \in L_2(\Omega)$
weakly on $\Omega$ and there exists a $\psi \in L_2(\widetilde \Gamma)$ 
such that $(B^* \psi - \alpha \, \Tr u)$ is the $(a_{ij})$-normal derivative of $u$
and $A j(u) = (\ca u, \psi)$.
Since $A j(u) \in j(H^1_{b,\sigma}(\Omega))$ one deduces that 
$\ca u \in H^1_{b,\sigma}(\Omega)$ and $j(\ca u) = A j(u) = (\ca u, \psi)$.
In particular, $\psi = B(\Tr \ca u)$.
Therefore $(B^* B (\Tr \ca u) - \alpha \, \Tr u)$ is the $(a_{ij})$-normal derivative of $u$
and $u \in D(A_1)$.
Then $A_1 u = \ca u = j^{-1} A_0 j(u)$.
Conversely, suppose that $u \in D(A_1)$.
Then $j(u) \in D(a_c)$ and it follows from Proposition~\ref{psesqui470}
that $j(u) \in D(A)$ with $A j(u) = (\ca u, B \Tr \ca u) = j(\ca u)$.
So $j(u) \in D(A_0)$ and $u \in D(j^{-1} A_0 j)$.\hfill$\Box$

\ruimte

In case of the Laplacian, i.e.\ if $a_{ij} = \delta_{ij}$,
the set $\Omega$ is bounded and Lipschitz, and if $B$ is the 
multiplication operator with a bounded measurable function $\beta$, 
then $D(A_1)$ is the set of all $u \in H^1(\Omega)$ such that 
$\Delta u \in H^1(\Omega)$ weakly on $\Omega$ and the 
normal derivative satisfies
\[
\frac{\partial u}{\partial \nu} + |\beta|^2 \, \Tr (\Delta u) + \alpha \, \Tr u = 0
.  \]
Moreover, $A_1 u = - \Delta u$.
Cf.\ \cite{AMPR} Remark~2.9.

\subsection*{Acknowledgement}

The authors thank Hendrik Vogt for critical comments
on earlier versions of Corollary~\ref{csesqui201.2} and 
Lemma~\ref{lsesqui410}.
In addition the authors thank the referee for his careful reading of the 
manuscript and his many improving comments.
The first named author is most grateful for the hospitality 
during a research stay at the University of Auckland.
Both authors highly appreciated the stimulating 
atmosphere and hospitality 
at the University of Torun in the framework of the 
Marie Curie `Transfer of Knowledge' programme TODEQ.

\end{document}